\pdfoutput=1

\documentclass[11pt]{article}
\pdfpagewidth=\paperwidth
\pdfpageheight=\paperheight
\usepackage[utf8]{inputenc}
\usepackage{url,geometry}
\usepackage[symbol]{footmisc}
\usepackage[final]{changes}
\usepackage{bbm}
\usepackage{bm}
\usepackage{amsmath, scalerel, amsthm, amssymb}
\usepackage{lscape}
\usepackage{algorithm}
\usepackage{adjustbox}
\usepackage{graphics}
\usepackage{graphicx}
\usepackage{wrapfig}
\usepackage{multirow}
\usepackage{pict2e,color}
\usepackage{amsfonts}
\usepackage{verbatim}
\usepackage[figuresright]{rotating}
\usepackage{booktabs}
\usepackage{url}
\usepackage{longtable}
\usepackage{tabularx}
\usepackage{multicol}
\usepackage{enumitem}
\usepackage{stmaryrd}
\usepackage[normalem]{ulem}
\usepackage{caption}
\usepackage{float}
\usepackage{tikz}
\usepackage{pgfplots}
\usepackage{setspace}
\usepackage{subfig}
\usepackage{algpseudocode}
\usepackage{xurl}
\usepgfplotslibrary{statistics}
\pgfplotsset{compat=1.17}

\usepackage[natbibapa]{apacite}
\renewcommand{\citet}[1]{\citeauthor{#1} \cite{#1}}

\newcommand{\Hline}{\hline}

\usetikzlibrary{pgfplots.groupplots}
\usetikzlibrary{shapes.misc, positioning}
\usetikzlibrary{calc,patterns,angles,quotes}
\usetikzlibrary{shapes,arrows,shapes.multipart}

\pgfplotsset{select coords between index/.style 2 args={
    x filter/.code={
        \ifnum\coordindex<#1\fi
        \ifnum\coordindex>#2\fi
    }
}}

\newcommand{\symfootnote}[1]{%
\let\oldthefootnote=\thefootnote%
\stepcounter{mpfootnote}%
\addtocounter{footnote}{-1}%
\renewcommand{\thefootnote}{\fnsymbol{mpfootnote}}%
\footnote{#1}%
\let\thefootnote=\oldthefootnote%
}

\renewcommand{\citet}[1]{\citeauthor{#1} \cite{#1}}

\newtheorem{theorem}{Theorem}
\newtheorem{lemma}[theorem]{Lemma}

\onehalfspace 

\title{The drop box location problem}
 
\author{\textbf{Adam Schmidt and Laura A.~Albert\footnote{Corresponding author}}\\
University of Wisconsin-Madison\\
Industrial and Systems Engineering \\
1513 University Avenue \\
Madison, Wisconsin 53706 \\
Email: \url{laura@engr.wisc.edu} 
}
\date{\today}

\renewcommand{\cite}{\citep}

\begin{document}

\maketitle

\begin{abstract}
For decades, voting-by-mail and the use of ballot drop boxes has substantially grown within the United States (U.S.), and in response, many U.S.~election officials have added new drop boxes to their voting infrastructure. However, existing guidance for locating drop boxes is limited. In this paper, we introduce an integer programming model, the drop box location problem (DBLP), to locate drop boxes. The DBLP considers criteria of cost, voter access, and risk. The cost of the drop box system is determined by  the fixed cost of adding drop boxes and the operational cost of a collection tour by a bipartisan team who regularly collects ballots from selected locations. The DBLP utilizes covering sets to ensure each voter is in close proximity to a drop box and incorporates a novel measure of access to measure the ability to use multiple voting pathways to vote. The DBLP is shown to be NP-Hard, and we introduce a heuristic to generate a large number of feasible solutions for policy makers to select from a posteriori. Using a real-world case study of Milwaukee, WI, U.S., we study the benefit of the DBLP. The results demonstrate that the proposed optimization model identifies drop box locations that perform well across multiple criteria. The results also demonstrate that the trade-off between cost, access, and risk is non-trivial, which supports the use of the proposed optimization-based approach to select drop box locations.
\end{abstract}

\renewcommand{\thefootnote}{\arabic{footnote}}
\section{Introduction}\label{sec:intro}

%
During the 2020 General election within the United States (U.S.), a record 46\% of U.S. voters cast a ballot by mail or absentee in-person \cite{mit_election_data__science_lab_voting_2021}.
Approximately 41\% of these voters cast a ballot using a drop box \cite{noauthor_sharp_2020}, which are temporary or permanent fixtures similar to United States Postal Service (USPS) postboxes.
Many states increased the number of drop boxes during 2020 in response to increased use of the vote-by-mail system and to help mitigate health risks associated with in-person voting \cite{corasaniti_postal_2020}. 
In total, forty states and Washington, D.C.~allowed some use of ballot drop boxes \cite{huord_where_2020}.
%
However, the increase in drop box use is likely not a one time event. 
The use of non-traditional voting methods within the United States has steadily grown since 1996 \cite{scherer_majority_2021}.
A recent survey of Wisconsin, U.S. election clerks found that approximately 78\% of election clerks would like some use of ballot drop boxes in future elections, and this percentage is higher among clerks from jurisdictions with a large voting age population \cite{burden_experiences_2021}.
Many states have since introduced legislation to expand the number of drop boxes available to voters\footnote{There are challenges to some proposals and even calls to restrict the use of these resources \cite{vasilogambros_lawmakers_2020}.} \cite{vasilogambros_lawmakers_2020}.

Reasons for casting a ballot using a drop box include the perceived security they offer, anticipated mail delays, and a lack of voter confidence in the USPS \cite{nw_2020_2020}.
For many individuals, drop boxes are also in close proximity of their home, work, or daily routine \cite{stewart_survey_2016}. 
Arguably, the primary benefit of drop boxes is the increased accessibility they offer to the voting infrastructure compared to in-person voting.
%
%
Studies suggest that
adding drop boxes to a voting system can increase voter turnout \cite{collingwood_drop_2018, mcguire_does_2020}.
\citet{mcguire_does_2020} found that a decrease in one mile to the nearest drop box increases the probability of voting by 0.64 percent. 
This finding aligns with the hypothesis of election participation first offered by \citet{downs_economic_1957}. 
According to this hypothesis, potential voters decide whether to vote by comparing the cost (e.g., time) of voting and the potential benefits from voting. 
It was later argued that voting cost is the significant driver of voter turnout \cite{sigelman_cost_1982,haspel_location_2005}.
We posit that the election infrastructure plays a large role in determining the cost to vote \cite{cantoni_precinct_2020, mcguire_does_2020, collingwood_drop_2018}.
Thus, if we can improve the accessibility of ballot drop boxes to voters by appropriately designing the drop box infrastructure, then we can increase voter participation, particularly among groups who previously had a high cost to vote and low turnout.

Although drop boxes can increase voter participation, 
there are many challenges associated with identifying drop box locations and managing the drop box voting system. 
First, drop boxes can pose a large financial cost.
Drop boxes can cost \$6,000 \cite{CISA_ballot_2020}, and designated video survillance cameras that increase drop box security can cost up to \$4,000 \cite{schaefer_ri_2020}.
%
Second, with an increased number of drop boxes, substantial time and resources must be devoted to collecting ballots.
During the election period, it is recommended that bipartisan teams regularly collect ballots \cite{CISA_ballot_2020}. 
If drop boxes are not strategically placed or if there are a large number of drop boxes, this route may be costly and leave less time to devote to other election tasks. 
Third, there are security risks associated with  ballot drop boxes that must be addressed \cite{scala_evaluating_2021},
although  drop boxes are considered reliable \cite{elections_project_staff_drop_2020,scala_evaluating_2021}. 
%
If the drop box specific security risks are mitigated appropriately, adding drop boxes to a voting system makes an adversarial attack on the electoral process more challenging.
This improves the overall security of the voting system, since the system becomes more distributed  \cite{scala_evaluating_2021}.
%
In addition to the previously mentioned challenges, elections are administered by state and local governments within the U.S., and each may have different voting processes.
While the vote-by-mail process is typically similar across difference jurisdictions within the U.S., each jurisdiction may have unique challenges or preferences that necessitates a detailed analysis of potential drop box system design.

In light of these complexities, existing guidelines for selecting drop box locations are often insufficient to support election administrators.
%
%
In 2020, the Cybersecurity and Infrastructure Security
Agency \cite{CISA_ballot_2020} recommended that a drop box be placed at the primary municipal building,
there be a drop box for every 15,000–20,000 registered voters, and more drop boxes should be added where there may be communities with historically low absentee ballot return rates.
However, these guidelines are not prescriptive enough to support administrators in  identifying an appropriate portfolio of drop box locations.
%
To our knowledge, the only analytical approach to selecting drop box locations
uses a Geographic Information System (GIS) to determine the locations that served the most voters, allowing for a maximum drive time of 10 minutes \cite{greene_vote-by-mail_2015}. 
This approach overlooks many of the trade-offs within the voting system and ignores socioeconomic differences between voters that may make voting more challenging for some individuals.

Without adequate decision support tools, election administrators may ultimately select drop box locations that perform poorly across multiple criteria by which voting systems are measured.
%
%
In this paper, we propose an integer program (IP) to support election administrators in determining how ballot drop boxes should be used in their voting systems when allowed by law\footnote{The ability to use or not use drop boxes and in what capacity is typically set by state law.}.
%
We formalize the IP as the  drop box location problem (DBLP).
To our knowledge, the DBLP is the first mathematical model of the ballot drop box system to support election planning.
The DBLP seeks to minimize the capital and operational cost of the drop box system,  ensure equity of access to the voting system, and mitigate risks associated with the drop box system. 
Loosely, we let \emph{access} refer to the proximity of the voting infrastructure (e.g., polling places, drop boxes) to voters and the ease with which voters can cast a ballot.
Expanding access through the use of drop boxes is an important aspect of the DBLP, since voter turnout is highly correlated with the distance needed to travel to cast a ballot \cite{cantoni_precinct_2020}.
We measure access to the drop box voting system using conventional covering sets.
In addition, we propose a function based on concepts from discrete choice theory to measure the level of access a voter has to the multiple voting pathways offered by the voting system.
The remainder of the paper is structured as follows.
In Section \ref{sec:lit}, we review the management science literature related to elections.
In Section \ref{sec:model}, we discuss measures by which the ballot drop box system can be assessed.
We then formalize the drop box location problem (DBLP) and introduce an 
IP formulation of the DBLP.
In Section \ref{sec:solnmethods}, we discuss solution methods for the DBLP and introduce a heuristic to quickly generate a collection of feasible solutions for election officials to select from a posteriori.
In Section \ref{sec:casestudy}, we introduce a case study of Milwaukee, WI, U.S. using real-world data. Using this case study, we demonstrate the value of our integer programming approach compared to rules-of-thumb that may otherwise be used.
We find that the DBLP outperforms the rules-of-thumb with respect to nearly all criteria considered.
We then investigate the trade-off between cost, access, and risk within potential drop box system designs.
We find that the trade-off is non-trivial, and the optimization-based approach provides value.
We conclude with a brief discussion in Section \ref{sec:discussion}.

\section{Literature Review}\label{sec:lit}

Much of the management science literature aimed at supporting election planning focuses primarily on in-person voting processes.
Some research focuses on identifying and describing the in-person voting process including quantifying the  arrival rate of in-person voters, the attrition rate of polling place queues, the check-in service rate, the time to vote, and poll worker characteristics  \cite{spencer2010long,stein_waiting_2020}. 
This research also studied how voting requirements (e.g., the introduction of voting identification requirements) impacts voting times \cite{stein_waiting_2020}.
Queueing theory has been widely used to analyze lines at polling locations and identify mitigating practices to avoid long lines \cite{stewart2015waiting,schmidt_designing_2021}. 
Since voting machines have been recognized as a bottleneck in the in-person voting process \cite{yang2009all}, a stream of papers has focused on the allocation of voting machines to polling locations \cite{allen_mitigating_2006,wang2015efficiency,edelstein2010queuing}.

Other research has focused on risks of voting systems rather than operational design.
The Election Assistance Commission (EAC) \cite{eac_elections_2009} analyzed threats to voting processes in the U.S. 
for seven voting technology types.
\citet{scala_evaluating_2021} identified security threats for mail-in voting processes and offered a relative score for each to identify the most important threats to address. 
They identify three drop-box related threats.
First, a misallocation of drop boxes can suppress voter turnout.
Second, a drop box can be damaged or destroyed.
Third, ballots within a drop box can be stolen or manipulated. 
They find the likelihood of drop box risks to be relatively low compared to other risks \cite{scala_evaluating_2021}. 
\citet{fitzsimmons_selecting_2020} study geographic-based risks by  introducing a control problem to study how voter turnout can be manipulated through the strategic selection of polling locations. 
%
A few papers attempt to detect disruptions or security incidents following an election \cite{highton_long_2006,allen_mitigating_2006}. 

There are no known papers intended to support election administrators in planning and managing the vote-by-mail system. 
%
Our proposed integer program addresses the risks of the drop box system \cite{scala_evaluating_2021} and employs concepts from the facility location literature.
Facility location problems are defined by a set of demands (e.g., voters) and a set of facilities (e.g., drop boxes) that can serve the needs of the demands.
Arguably the most widely used facility location model is the maximal covering location problem (MCLP)  \cite{church_maximal_1974}. 
In the MCLP, a demand is ``covered'' by, or can be served by, a predetermined set of locations called the \emph{covering set}. Facility locations are selected to maximize the number of demands covered by at least one facility.
The location set covering problem (LSCP) instead requires that all demands are covered and the cost of the selected facility locations is minimized \cite{toregas_location_1971}. 

The IP introduced within this paper extends the covering tour problem (CTP) \cite{gendreau_covering_1997}, which is a variant of the LSCP, by considering additional constraints and objective function terms.
These changes allow us to accurately model the drop box voting system.
A CTP instance is defined by an undirected weighted graph with two mutually exclusive and exhaustive set of nodes, the \emph{tour} nodes and \emph{coverage nodes}.
The objective of the CTP is to find a Hamiltonian tour of minimal length over a subset of the tour nodes such that each coverage node is covered by at least one node visited by the tour. 
The CTP is NP-Hard since the traveling salesman problem (TSP) can be reduced to it \cite{gendreau_covering_1997}.
Several solution methods, including exact \cite{gendreau_covering_1997,baldacci_scatter_2005} and heuristic \cite{murakami_generalized_2018,vargas_dynamic_2017}, have been proposed for the CTP. 
This paper represents the first known application of a CTP variation to voting systems.

\section{Problem Definition}\label{sec:model}
\newcommand\drawredrect{%
  \begin{tikzpicture}                                       
  \tikz\draw[red,thick,dotted] (0,0) rectangle (2ex,1.25ex);      
  \end{tikzpicture}%
}  

Election administrators in the U.S.~face many questions regarding the use of ballot drop boxes including whether drop boxes should be added to their local voting system and how drop boxes may affect voting performance measures.
If election administrators decide to add drop boxes, they must decide how many drop boxes to add and where they should be located.
The DBLP introduced in this section identifies the optimal placement of drop boxes once election administrators decide to add drop boxes to the voting system.
However, election administrators can use the model during the election planning process to assess the cost, access, and risk of a potential drop box system. This can inform their decision of whether or not to add any drop boxes to the voting system.

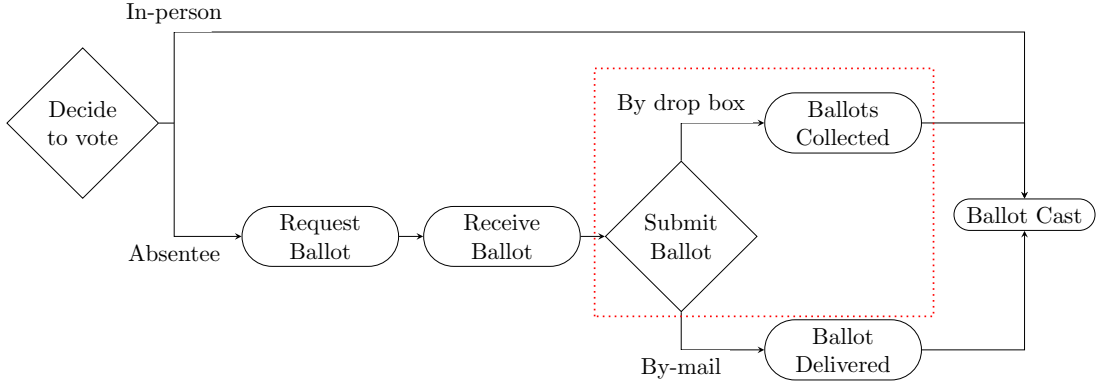
\begin{figure}[htb!]
    \centering
    \begin{adjustbox}{width=.9\textwidth}
    \begin{tikzpicture}
    
    \node[draw, diamond, text width=1.25cm, align = center ] (1) at (0,0) {Decide to vote};
        
    \begin{scope}[every node/.style={draw, rounded rectangle, text width=2cm, align = center  }]
        \node[below right=0.75cm and 2.5cm of 1] (3) {Request Ballot};
        \node[right=0.4cm of 3] (4) {Receive Ballot};
    \end{scope}   
    
     \node[draw, diamond, text width=1.25cm, align = center , right=0.4cm of 4] (5) {Submit Ballot};
       
    \begin{scope}[every node/.style={draw, rounded rectangle, text width=2cm, align = center  }]   
        \node[above right=0.75cm and 1.25cm of 5] (6) {Ballots Collected};
        \node[below right=0.75 and 1.25cm of 5] (7) {Ballot Delivered};
        \node[below right=0.75 and 1.3 of 6] (8) {Ballot Cast};
    \end{scope}
    
    \node[right=0.25cm of 1] (0) {};
    
    \node[above right=0.75cm and 1.5cm of 1] (2) {};
    
    \begin{scope}[>={stealth[black]},
                  every node/.style={fill=white},
                  every edge/.style={draw=black}]
       \draw [-] (1.east)  edge  (0);
       \draw [-] (0.west)  |- node[above]{In-person} (2);
       \draw [->] (0.west)  |- node[below]{Absentee} (3);
       \path [->] (3) edge (4);
       \path [->] (4) edge (5);
       \draw [->] (5)  |- node[above] {By drop box} (6);
       \draw [->] (5)  |- node[below] {By-mail} (7);
       \draw [->] (6)  -| (8);
       \draw [->] (7)  -| (8);
       \draw [->] ($(2.west)$)  -| (8);
    \end{scope}
    
    \draw[red,thick,dotted] ($(5.south west)+(-0.8,-0.7)$)  rectangle ($(6.north east)+(0.7,+0.4)$);
\end{tikzpicture}
\end{adjustbox}
    \caption{Typical pathways to cast a ballot, divided into in-person and absentee, and the component corresponding to the use of ballot drop boxes (\protect\drawredrect).}
    \label{fig:voting}
\end{figure}

The decisions surrounding the use of drop boxes are complex due to the number of potential locations for drop boxes, concerns about equity within the voting process, and the multiple criteria by which voting systems are measured.
The most widely reported election performance metrics in the U.S.~are the number of individuals registered to vote and the fraction of eligible voters that cast a ballot, known as \emph{voter turnout} \cite{mit_election_data__science_lab_elections_2022}. 
In most states, there are multiple pathways by which voters can cast a ballot, and the accessibility of each pathway can influence voter turnout.
Figure \ref{fig:voting} describes the two main pathways, which are typically divided into `in-person' or `absentee'.
With in-person voting, a voter
obtains and casts a ballot at their assigned polling location, typically on election day.
With absentee voting\footnote{Some states, such as Washington, use the ``absentee'' voting process as their primary voting method. Thus, we use ``absentee'' loosely in this paper, and sometimes refer to it as the vote-by-mail process.}, a voter requests a ballot be sent to them and the completed ballot is then returned either through the mail or using a drop box.
In some states, voters must provide a reason to vote absentee, while in 34 states there is ``no-excuse'' absentee voting \cite{national_conference_of_state_legislatures_voting_2022}.

In addition to voter-based election metrics, the cost and security of the voting system is a key concern. 
The cost of an election is comprised of both infrastructure-based costs (e.g., polling locations) and resource based costs (e.g., staff).
The security of a voting system is not typically measured or reported to the public, despite being a major concern of officials and the public.

In this paper, we are concerned with a sub-pathway of the vote-by-mail process where the voter submits a ballot using a drop box.
In this pathway, a voter first requests and receives a ballot through the mail.
They then decide to submit a ballot using a drop box rather than through the mail (or not returning it at all).
This decision is influenced by the proximity of a drop box to the voter and the distrust the voter has in the USPS \cite{mcguire_does_2020}.
A team of poll workers then collects ballots from the drop boxes, and the ballots are processed at an official election building.
In this paper, we focus on the system related to the steps outlined in red  (\protect\drawredrect), since they are the steps that are unique to the drop box system and are influenced by the locations of the drop boxes.


\subsection{Assessing Drop Box Infrastructure}\label{sec:metrics}
There are two metrics typically used to assess the vote-by-mail system: the proportion of requested ballots that are returned and the number of ballots rejected \cite{mit_election_data__science_lab_elections_2022}.
The use of ballot drop boxes can lower the rejection rate of mail ballots by reducing the time it takes a ballot to return to election officials. As a result, a voter can be notified of an incorrectly marked ballot more quickly to allow the voter to resubmit their ballot before the election deadline.
This is a benefit that we do not explicitly consider in our model.
We also posit based on empirical research that a well-designed drop box system can lead to a higher proportion of returned mail ballots and a higher voter turnout by improving the accessibility of the voting infrastructure \cite{downs_economic_1957,mcguire_does_2020}.

We elaborate on how access to the voting system is measured.
We employ the concept of \emph{coverage} to measure the access voters have to the drop box system.
Under the concept of coverage, a voter covered by a selected drop box location is assumed to have access to the drop box voting system.
The locations that provide a voter coverage are called its \emph{covering set}.
Covering sets are flexible and can be defined to account for different modes of transportation, vehicle ownership, and other socioeconomic factors.
However,
drop boxes are a subcomponent of a larger voting system, and coverage overlooks the access provided by non-drop box voting pathways.
In reality, some individuals may have better access to in-person voting than others, and adding drop boxes near them may not substantially benefit them.
This necessitates a second measure of access that distinguishes access to the complete voting infrastructure from coverage by the drop box system.

We introduce an access function based on the multinomial/conditional logit model from discrete choice theory \cite{aloulou_application_2018} to capture this dynamic. 
The application of discrete choice theory to questions within political science is most commonly used to explain or predict choices within a multi-candidate (or party) election \cite{glasgow_discrete_2008}.
Discrete choice models have also been used to predict how individuals interact with infrastructure in different application domains. One of the earliest cases of this was the application of a conditional logit model to predict the use of the Bay Area Rapid Transit prior to its construction \cite{train_discrete_2009}.
To the best of our knowledge, our paper represents the first use of a function based on discrete choice theory to model access within an optimization model.

The function we introduce makes use of some parameters.
Let $v^1_{w} > 0$ be a measure of accessibility\footnote{In its exact form, $v^1_{w} = e^{U_w^1}$ where $U_w^1$ represents the utility of voting using the non-drop box voting system.} to the non-drop box voting system (e.g., in-person polling locations) for voters $w$. This can be determined, for example, by the distance to the nearest polling location. 
Let $a_{nw} > 0$ be a measure of the access\footnote{In its exact form, $a_{nw} = e^{U_{wn}}$ where $U_{wn}$ represents the utility of voting by using drop box $n$.} that a drop box at location $n$ would provide to $w$. 
This can be determined in part by the proximity of the location to the voters' places of residence and work and by the various transportation modes available between the voters and the drop box location.
Based on empirical studies, the value of $a_{nw}$ should be increasing with decreasing distance \cite{mcguire_does_2020}.
Finally, let $v^0_{w} > 0$ be the propensity of $w$ not to vote\footnote{In its exact form, $v^0_{w} = e^{U_w^0}$ where $U_w^0$ represents the utility of not voting.}.
This could be informed by the historical non-voting rate (complement of turnout) or using surveys.
Using these parameters, we introduce the following \emph{access function} to measure the access a group of individuals $w$ has to all voting pathways where $N^*$ represents the set of selected drop box locations:
\begin{align}
        A_w(N^*) : = \frac{v^1_{w}+\sum_{j \in N^*} a_{jw} }{v^0_{w}+v^1_{w}+\sum_{j \in N^*} a_{jw} } \nonumber 
\end{align}
The access function takes values between zero and one. A value closer to one means that the voting system, including the new ballot drop boxes, is more accessible to individuals $w$, whereas value closer to zero means that the voting system is relatively inaccessible to individuals $w$.
In this way, a higher access function value suggests higher turnout for $w$.

The access function can still be used when a strict interpretation is not reasonable or is not feasible due to data availability, since the benefit of the access function is a result of its structure.
First, the access function models access as a non-binary measure.
Second, adding any drop box to the voting system increases the value of the access function but to varying degrees based on the locations of the voter and the drop box.
Third, each voter has some heterogeneous level of access to non-drop box voting methods captured by $v_w^1$, and this access is treated as a constant within the scope of the decision to location drop boxes.
Each voter also has a heterogeneous access function value when no drop boxes are added to the voting system, $A_w(\emptyset) = \frac{v_w^1}{v_w^0 + v_w^1}$, which is reflective of heterogeneous turnout rates.
Fourth, the benefit of adding a drop box near a voter is marginally decreasing as 
the access function value increases.
This incentivises placing drop boxes near populations with low levels of access to other voting pathways.

While it is desirable to increase voter turnout and access to the voting system, expanding the use of ballot drop boxes may increase the financial cost of managing the election.
The costs of the ballot drop box system can be broken into two major groups: fixed or operational. 
\emph{Fixed costs} represent the ``per drop box'' costs such as the initial purchase and costs of securing and maintaining the drop box.
Each location may have a different fixed cost due to varying installation and security equipment requirements.
Once drop boxes are installed, jurisdictions incur an \emph{operational cost} for a bipartisan team to collect ballots from the drop boxes \cite{CISA_ballot_2020}.
The operational cost is determined, in part, by the distance between drop boxes, the opportunity cost of bipartisan team's time, and the frequency at which the ballots are collected during an election.
We assume that a bipartisan team collects ballots from all drop boxes whenever a collection is conducted, and the drop boxes are visited in an order that minimizes the operational cost, referred to as the \emph{collection tour}.

In addition to introducing new financial costs, drop boxes introduce three types of risks to the voting process that can be mitigated through design requirements.
The first risk is that ballot drop boxes can be misallocated in a way that causes voter suppression \cite{scala_evaluating_2021}.
There are two components to this risk.
The first is the potential to misallocate drop boxes such that access to the drop box voting system is inequitable.
The second is the potential to misallocate drop boxes such that the access to the entire voting system defined by the multiple voting pathways is inequitable.
These risks are reflected by the number of voters covered by a drop box, using the same definition of coverage introduced earlier,  and the value of the access function for each voter, respectively.
We can mitigate the risk of voter suppression by requiring that each voter is covered by at least one drop box and that the value of the access function meets some minimal threshold for all voters.

The second risk is that a drop box could be damaged or destroyed \cite{scala_evaluating_2021}.
A nefarious actor could influence an election by targeting drop boxes that provide access to certain voters.
The impact of this risk can be mitigated by requiring all voters to be covered by multiple drop boxes, so that voters have redundant access to the drop box system.
%

The last risk is that ballots submitted to a drop box could be stolen or manipulated.
The impact of this risk can be mitigated by ensuring 
that the collection tour has a low cost.
When the collection tour has a low cost,  election officials can collect ballots often, leaving fewer at risk.
Other implicit design choices also mitigate this third risk.
For example, requiring a bipartisan team to collect ballots, rather than one individual, reduces the risk of an insider attack.
Likewise, incorporating security related costs, such as the cost of a video surveillance system, into the fixed cost of a drop box mitigates the risks associated with it.

There are additional risks and mitigations associated with the voting process that are not unique to the drop box infrastructure. For example, there is a risk of an insider attack on ballots stored at an election building after collected from the drop boxes \cite{scala_evaluating_2021}.
However, these additional risks are outside the scope of the system  considered in this paper (see Figure \ref{fig:voting}).

\subsection{The  Drop Box Location Problem (DBLP)}\label{sec:DBLP}

We now formally introduce an IP formulation of the drop box location problem (DBLP) using the sets, parameters, and variables presented in Table \ref{drop box:T:notation}.
%

\begin{table}[hbt!]
\centering
\caption{Notation}\label{drop box:T:notation}
\begin{tabular*}{\columnwidth}{@{}l@{\extracolsep{\fill}}c@{\extracolsep{\fill}}p{0.8\columnwidth}}
\multicolumn{3}{l}{\textbf{Sets}}\\
\Hline
$W$ & = & the set of voter populations \\
$N$ & = & the potential drop box locations \\ 
  $T \subseteq N$ & = & the locations at which a drop box must be placed \\ 
 $E$ & = & all pairs $i \in N$, $j\in N$ such that $i \neq j$ and $(j,i)\neq E$ \\
$N_w \subseteq N$ & = &drop box location that cover $w \in W$, $|N_w| \geq 2$ \\
\Hline
\\
\end{tabular*}
\begin{tabular*}{\columnwidth}{@{}l@{\extracolsep{\fill}}c@{\extracolsep{\fill}}p{0.8\columnwidth}}
\multicolumn{3}{l}{\textbf{Parameters}}\\
\Hline
$s$ & = & the start and end of the collection tour \\
$f_j$ & = & the fixed cost of placing a drop box at location $j \in N$ \\
$c_{ij}$ & = & the operational cost of traveling between $i \in N$ and $j \in N$ in the collection tour \\
$v^0_{w}$ & = & the propensity of $w \in W$ not to vote \\
$v^1_{w}$ & = & the accessibility of the non-drop box voting system to $w \in W$ \\
$a_{jw}$ & = & the accessibility of location $j \in N$ to $w \in W$ \\
$r$ & = & minimal allowable value for the access function \\
$q$ & = & minimal number of drop boxes covering each $ w \in W$ \\
\Hline
\\
\end{tabular*}
\begin{tabular*}{\columnwidth}{@{}l@{\extracolsep{\fill}}c@{\extracolsep{\fill}}p{0.8\columnwidth}}
\multicolumn{3}{l}{\textbf{Decision Variables}}\\
\Hline
$x_{ij}$ & = & 1 if the collection tour moves between $i$ and $j$ $(i,j) \in E$ and 0 otherwise \\
$y_j$ & $=$ & 1 if a drop box is placed at location $j \in N$ and 0 otherwise  \\
\Hline
\\
\end{tabular*}\vspace{-0.5cm}
\end{table}

%
The DBLP selects drop box locations from a set of potential locations, $N$. Potential drop box locations can be identified using existing guidelines \cite{CISA_ballot_2020, mcguire_does_2020}. 
Let $y_n$ be a decision variable that equals one
if a drop box is located at location $n \in N$ and zero otherwise.
Once drop box locations are selected, a collection tour over them must be found to determine the operational cost of the drop box system.
Let $x_{ij}$ be a decision variable
that equals one if the collection tour travels between drop box $i$ and drop box $j$, $(i,j) \in E$, and zero otherwise, where $E$ represents  all pairs $i \in N$, $j\in N$ such that $i \neq j$ and $(j,i)\neq E$ .
We assume the collection tour always begins and ends at a drop box\footnote{When a drop box is not located at $s$, then the model is still valid. Simply let $f_s = 0$ and $a_{sw} = 0$ for each $w \in W$, while $s$ is not a member of any covering set.} located at $s$ (e.g., primary municipal building).
Let $T$
represent the locations at which there must be a drop box within our solution (e.g., existing drop box locations).
The set $T$ is always non-empty, since $T=\{s\}$ in the extreme case.
%
For each location $j \in N$, let $f_j$ equal the fixed cost of 
a drop box at $j$.
Let $c_{ij}$ represent the operational cost of traveling between drop boxes $(i,j) \in E$ on the collection tour.

%
Using this notation, we formalize the three goals of the DBLP.
The first goal is to minimize the total cost associated with the selected drop box locations.
The total cost of the drop box system is the sum of the fixed costs and the cost of the collection tour,  $z_1 := \sum_{j \in N} f_j y_j + \sum_{(i,j)\in E} c_{ij} x_{ij} $.
The value of $z_1$ serves as the objective\footnote{Election administers likely have a fixed budget, but the amount allocated to managing the drop box system is likely not predetermined. Thus, we wish to minimize the proportion of the budget allocated to the drop box system.} function in the IP formulation of the DBLP.


The second goal of the DBLP is to equitably improve the accessibility of the voting system.
Let $W$ denote the collection of voter populations.
Let $N_w \subseteq N$ represent the drop boxes that cover $w \in W$.
We ensure equitable access to the drop box system\footnote{When $q \geq 1$.} by requiring that each voter is covered by $q$ drop boxes.
Reasonable values of $q$ are $0$, $1$, or $2$.
The cardinality of each covering set must be at least $q$, $|N_w| \geq q$ for all $w\in W$, otherwise the problem is infeasible.
We ensure equitable access to all voting pathways by requiring that the access function value is at least $r$ for each $w\in W$, $\min_{w \in W} A_w (N') \geq r$ where $N' = \{n \in N : y_n = 1\}$ are the selected drop box locations.
This constraint can be viewed as a second objective for the DBLP using the epsilon-constraint approach for multi-objective optimization problems \cite{mavrotas_effective_2009}.


The third goal of the DBLP is to mitigate the risks associated with the drop box voting system. 
The risk of misallocating drop boxes in a way that leads to voter suppression is addressed by the second goal of the DBLP.
The risk of ballots being susceptible to manipulation once submitted to a drop box is addressed by minimizing the cost of the collection tour, which is captured within $z_1$.
The risk of damage to or destruction of drop boxes is a way that degradates voter access to the voting system is mitigated by ensuring each voter is covered by $q$ drop boxes when $q \geq 2$.

If the optimal solution to the DBLP locates two or fewer drop boxes\footnote{It can be easily checked whether two or fewer drop boxes are needed to satisfy the constraints of the model.}, the collection tour visiting the drop box locations is trivial. 
Thus, we assume that at least three drop boxes are be selected in the optimal solution. Under this assumption, we can formulate the DBLP using the following IP.
\begin{align}
    \underset{x,y}{\min} \ &  z_1 =  \sum_{(i,j) \in E} c_{ij} x_{ij} + \sum_{j \in N} f_j y_j \label{model:obj1} \\
    \text{s.t.} \ 
    & r(v^0_{w}+v^1_{w}+\sum_{j \in N} a_{jw}y_{j} ) \leq v^1_{w}+\sum_{j \in N} a_{jw}y_{j} &  \forall \ w\in W  \label{model:obj2} \\
    & \sum_{j \in N_w} y_j \geq q & \forall \ w\in W \label{model:basecoverage}\\
    & y_j = 1 & \forall \ j \in T \label{model:existing} \\
    & \sum_{i \in N : (i,j) \in E} x_{ij} = 2y_j & \forall \ j \in N \label{model:balance} \\
    & \sum_{\substack{(i,j)\in E : i \in S, j \in N \setminus S \\ \hspace{0.55cm}\text{ or } j \in S, i \in N \setminus S}}    x_{ij} \geq 2y_t  & \substack{\forall S \subset N, \ 2 \leq |S| \leq |N|-2,\\ T\setminus S \neq \emptyset, t \in S} \label{model:subtourelim} \\
    &  y_{j} \in \{0,1\} & \forall \ j \in N \label{model:ybin}\\
    & x_{ij} \in \{0,1\} & \forall \ (i, j) \in E \label{model:xbin} 
\end{align}

The objective \eqref{model:obj1} is to minimize the total cost of the drop box system.
Constraint set \eqref{model:obj2} requires that the value of the access function is at least $r$ for each $w \in W$.
Constraint set \eqref{model:basecoverage} ensures that each $w \in W$  is covered by at least $q$ drop boxes within their respective covering set. 
Constraint set \eqref{model:existing} ensures that a drop box is added at each location in $T$.
Constraint sets \eqref{model:balance} and \eqref{model:subtourelim}  are used to determine the collection tour over the selected drop box locations using constraints originally introduced for the CTP \cite{gendreau_covering_1997}.
Constraint set \eqref{model:balance} ensures that each selected drop box location is visited by the collection tour exactly once. 
Constraint set \eqref{model:subtourelim} introduces subtour elimination constraints. 
Note that these constraints differ from the subtour elimination constraints commonly seen in the TSP, since not all locations $N$ must be visited by collection tour.
The bound on the summation refers to the edges in $E$ such that the edge is incident to one node in $S$ and one in $N\setminus S$.
Constraint sets \eqref{model:ybin} and \eqref{model:xbin} require the decision variables to be binary.

\subsection{Model Properties}

The DBLP is challenging to solve using existing solution techniques.
This idea is formalized in Theorem \ref{thm:nphard}, which states that the DBLP is NP-Hard. A proof is provided in the Supplementary materials \ref{appx:proofs}.

\begin{theorem}\label{thm:nphard}
The DBLP is NP-Hard.
\end{theorem}

In some instances, the DBLP integer program may be large due to a large number of voter populations.
We present a condition that allow us to determine when certain constraints from constraint set \eqref{model:obj2} can be removed from the IP, which reduces the size of the integer program instance and potentially reduces the time needed to find an optimal solution.
Lemma \ref{prop:covdom} gives a sufficient condition for which the constraint corresponding to a voter population $w \in W$ in constraint set \eqref{model:obj2} can be removed from the DBLP integer program, since the access function value is guaranteed to be smaller for another voter population $\hat{w} \in W$ for all choices of drop box locations.
 A proof is provided in the Supplementary materials \ref{appx:proofs}.



\begin{lemma}\label{prop:covdom}
Let $w,\hat{w} \in W$ be two voter populations. 
If the access function parameters satisfy $v^0_{\hat{w}} \geq v^0_{w}$, $v^1_{\hat{w}} + \sum_{n \in T} a_{n \hat{w}} \leq v^1_{w} + \sum_{n \in T} a_{n w} $, and $a_{n \hat{w}} \leq a_{n w}$ for each $n \in N\setminus T$, then for any subset of drop box locations $N' \subseteq N$, such that $T \subseteq N'$:
$$\frac{v^1_{\hat{w}}+\sum_{n \in N'} a_{n\hat{w}} }{v^0_{\hat{w}}+v^1_{\hat{w}}+\sum_{n \in N'} a_{n\hat{w}} } \leq \frac{v^1_{w}+\sum_{n \in N'} a_{nw} }{v^0_{w}+v^1_{w}+\sum_{n \in N'} a_{nw} }$$
\end{lemma}



This property may be satisfied in realistic instances of the DBLP.
Consider population $\hat{w}$ that lies on the exterior boundary of the jurisdiction. 
Consider a $w$ that lies just inside of $\hat{w}$ within the jurisdiction such that $w$ is closer than $\hat{w}$ to all potential drop box locations and polling locations. 
The voters in $w$ have higher access to the voting infrastructure than the voters in $\hat{w}$.
In this case, the properties of Lemmas \ref{prop:covdom} are likely to be satisfied.

\subsection{Model Variations}\label{sec:variations}
The DBLP is designed to determine drop box locations that satisfy the requirements for a drop box system in most jurisdictions.
However, election administrators may wish to explore solutions that are not identified by the standard formulation of the DBLP or may wish to tailor the model to their situation.
In this subsection, we discuss five modifications that can be made to the DBLP. 


First, the value of $q$ determines the number of drop boxes that must cover each voter, and determines the access that the voters have to the drop box system to an extent.
Letting $q = 0$, voters are not required to be covered by drop boxes. Instead, a cost effective set of drop box locations are selected such that all voters have a minimum level of access ($r$) to the voting system.
Letting $q = 1$, drop box locations are selected so that each voter is guaranteed access to the drop box system in addition to meeting a minimal level of access to all voting pathways  ($r$).
Letting $q = 2$, voters are guaranteed access to both the drop box and complete voting system in a way that also mitigates the impact the destruction of a drop box could have on voter access.

The second variation we consider is a change to the covering sets $N_w$ for $w \in W$.
The covering sets $N_w$ typically include locations within a predefined time or distance threshold from a voter.
Decreasing or increasing the time threshold used can make constraint set  \eqref{model:basecoverage} more or less restrictive, respectively.
When the covering sets are determined using a shorter time threshold, drop boxes within a voter's covering set are required to be located closer to the voter.
This may make the drop boxes more accessible to all voters, but also increases costs.
When the covering sets are determined using a larger time threshold, drop boxes are allowed to be located further away while still satisfying constraint set  \eqref{model:basecoverage}, which results in lower cost. 

Third,  we can replace the cost objective of the DBLP with other goals.
We can instead maximize the minimum access function\footnote{It is fairly straightforward to convert constraint set \eqref{model:obj2} into a linear equivalent by using one minus the access function value, which is an equivalent measure of access.} or maximize the number of voters covered by at least $q$ drop boxes.
In the latter case, we introduce a new indicator variable $\delta_w$ and the following constraint
\begin{align}
    q \delta_w \leq \sum_{j \in N_w} y_j \quad \forall \ w\in W
\end{align}
The objective is then to maximize $\sum_{w\in W} p_w \delta_w$ where $p_w$ represents the number of voters in $w$.
With this objective, we can use constraint set \eqref{model:basecoverage} to ensure each population $w \in W$ is covered by at least some $q'$ ($0 \leq q' < q$) drop boxes.
When $z_1$ is no longer the objective of the integer program, a constraint can be added to ensure that the cost of the drop box system is no more than some budget $B$,
\begin{align}
    \sum_{(i,j)\in E} c_{ij}x_{ij} + \sum_{n \in N} f_n y_n \leq B
\end{align}
With this constraint, feasibility of the DBLP is no longer guaranteed. Infeasibility can be informative to election administrators.
%


Fourth, election administrators may wish to restrict the cost of the collection tour, since they may have a limited  operational budget for collecting ballots (e.g., limited staff).
We can limit the cost of the collection tour to no more than $c_{\max}$ by introducing the following constraint
\begin{align}
    \sum_{(i,j)\in E} c_{ij} x_{ij} \leq c_{\max}
\end{align}

Fifth, election administrators may wish to locate a specific number, $k$, of drop boxes.
This may occur when they have already purchased drop boxes or they wish to add a drop box for every 15,000-20,000 registered voters as recommended by \cite{CISA_ballot_2020}.
This can be enforced by adding the following constraint
\begin{align}
    \sum_{n \in N} y_n = k
\end{align}

\section{Solution Methods}\label{sec:solnmethods}


Within this section, we present solutions methods for the original DBLP formulation.

\subsection{Objective Reformulation}\label{sec:obj1reform}

Constraint sets \eqref{model:basecoverage}-\eqref{model:subtourelim}
are similar to constraints that may be found in an integer program for the CTP \cite{gendreau_covering_1997}.
However, the objective of the CTP only considers operational costs \cite{gendreau_covering_1997}.
Thus, it is desirable to reformulate objective $z_1$ to preserve properties of the CTP within the DBLP.
We can then use components of solution methods  for the CTP within solution methods for the DBLP.

We present a reformulation of $z_1$ to remove the use of $y_n$ variables. 
Note that constraints (6) enforce that for any drop box $n$ visited by a feasible tour, there must be exactly two drop boxes visited before or after $n$. Thus, we can reformulate $z_1$ as follows:
\begin{align}
    z_1 & = \sum_{i,j \in E} c_{ij} x_{ij} + \sum_{j \in N} f_j y_j \nonumber \\ 
    & = \sum_{i,j \in E} c_{ij} x_{ij} +  \sum_{j \in N}\sum_{i \in N : i,j \in E} f_j  x_{ij}/2 \nonumber \\
    & =\sum_{i,j \in E} (c_{ij}+f_i/2+f_j/2)  x_{ij} \nonumber \\ 
    & =  \sum_{i,j \in E}\hat{c}_{ij}  x_{ij} \nonumber
\end{align}
where $\hat{ c}_{ij} := c_{ij}+f_i/2+f_j/2$ for each $(i,j) \in E$.
With this reformulation, $z_1$ takes the same form as the standard objective for the CTP.

\subsection{Lazy Constraint Method}\label{sec:exact}

Branch and bound is one of the most common techniques used to solve IPs, and we employ it to solve the DBLP.
However, constraint set \eqref{model:subtourelim} defines an exponential number of constraints, so we introduce a lazy constraint approach to solve the DBLP.
First, we solve the DBLP without constraint set \eqref{model:subtourelim}.
Once an optimal solution is found, we determine if any of the constraints from constraint set \eqref{model:subtourelim} are violated.
If so, we add in at least one violated constraint and resolve the IP using branch and bound.
Most modern optimization packages support the implementation of lazy constraints.
The reformulation of the objective introduced in Section \ref{sec:obj1reform} can be used throughout the procedure, but it is not required.

We introduce a new polynomial time algorithm, Algorithm \ref{alg:lazy}, to find violated inequalities from constraint set \eqref{model:subtourelim} given an $x^* \in \{0,1\}^{|E|}$.
The approach we take is adapted from an approach used for the TSP \cite{gurobi_optimization_tsppy_nodate} to account for the fact that not all potential drop box locations must be visited by the tour in the DBLP.
Algorithm \ref{alg:lazy} first finds all subtours defined by $x^*$ (line 1).  
Each subtour that does not include all required locations $T$ (lines 2-4) must be associated with least one violated constraint. 
For all\footnote{There is a trade-off between adding a constraint for all $t \in \hat{S}$ (increasing the size of the IP) and adding a constraint for a small number of elements in $\hat{S}$ (increasing the number of times the search procedure occurs).} locations $t$ visited by the subtour, we add the violated constraint (line 7).

\begin{algorithm}
  \caption{Lazy($x^*$)}\label{alg:lazy}
  \begin{algorithmic}[1]
    \State $H = $ collection of subtours defined by $x^*$
    \For{each subtour $h \in H$}
        \State $\hat{S} = $ drop box locations visited by $h$
          \If{$T \setminus \hat{S} \neq \emptyset$}
            \State \Return  $\sum_{(i,j)\in E : i \in \hat{S}, j \in N \setminus \hat{S} \text{ or } j \in \hat{S}, i \in N \setminus \hat{S}} x_{ij} \geq 2y_t$ for each node $t \in \hat{S}$ 
          \EndIf
    \EndFor
  \end{algorithmic}
\end{algorithm}

We comment on the correctness of Algorithm \ref{alg:lazy}.
Specifically, given an integer $x^* \in \{0,1\}^{|E|}$, Algorithm \ref{alg:lazy} finds a violated constraint from constraint set \eqref{model:subtourelim}, if one exists.
If a constraint is violated, there must exist a $S$ such that $S \subset N, \ 2 \leq |S| \leq |N|-2, \ T\setminus S \neq \emptyset$ and for some ${t^*} \in S$,
${ \sum_{(i,j)\in E : i \in \hat{S}, j \in N \setminus \hat{S} \text{ or } j \in \hat{S}, i \in N \setminus \hat{S}}    x^*_{ij} < 2y_{t^*}}$.
Since the left hand side of the inequality is at least zero, $t^*$ must represent a selected drop box location ($y_{t^*} = 1$).
Moreover, the feasibility of $x^*$ with regards to constraint set \eqref{model:balance} implies that $\sum_{(i,j)\in E : i \in \hat{S}, j \in N \setminus \hat{S} \text{ or } j \in \hat{S}, i \in N \setminus \hat{S}} x^*_{ij} = 0$.
Thus, ${t^*}$ must be a member of some subtour visiting locations $\hat{S} \subseteq S$.
The set $\hat{S}$ must contain at least three elements and can contain not more than $|N|-3$ elements as a result of constraint set \eqref{model:balance}.
Since $T\setminus S \neq \emptyset$, it is also true that $T\setminus \hat{S} \neq \emptyset$.
Thus, the existence of a $S$ implies the existence of a $\hat{S}$ whose elements form a subtour in $x^*$ such that $\hat{S} \subset N, \ 2 \leq |\hat{S}| \leq |N|-2, \ T\setminus \hat{S} \neq \emptyset$ and $ \sum_{(i,j)\in E : i \in \hat{S}, j \in N \setminus \hat{S} \text{ or } j \in \hat{S}, i \in N \setminus \hat{S}}    x^*_{ij} < 2y_{t}$ for all ${t} \in \hat{S}$.
Algorithm \ref{alg:lazy} identifies $ \hat{S}$ and returns the corresponding constraint.

\subsection{A Heuristic Method}\label{sec:heuristic}
The lazy constraint method can be used to find solutions to moderately sized problem instances, but large instances require long computational time.
Moreover, if an appropriate value for $r$ is not known by election administrators, the DBLP must be solved repeatedly to allow election administrators to select among possible solutions a posteriori\footnote{This may be desirable even if $r$ is believed to be known.}, which substantially increases the necessary computational time.
In this section, we present a heuristic that identifies multiple solutions to the DBLP, each corresponding to a unique value of $r$. 
Depending on the implementation, the heuristic is polynomially solvable.
We provide pseudocode for this heuristic in the Supplementary materials \ref{appx:heur}.
The heuristic requires the objective reformulation discussed in Section \ref{sec:obj1reform}.

The heuristic first identifies a feasible solution to the DBLP corresponding to $r = 0$.
When $q = 0$, we find a tour visiting the nodes of $T$ using any method for the TSP\footnote{Since we do not specify the methods to solve the TSP or CTP, this heuristic is in fact a \emph{heuristic scheme}.}.
When $q = 1$, the DBLP with $r = 0$ is equivalent to the CTP when the DBLP objective is reformulated as introduced in Section \ref{sec:obj1reform}. 
Any solution method for the CTP can be used to identify a feasible solution.
When $q = 2$, the DBLP with $r = 0$ is equivalent to the CTP when the DBLP objective is reformulated as introduced in Section \ref{sec:obj1reform}, except the CTP only requires single coverage of each $w \in W$.
We construct a solution that satisfies double coverage using any exact or heuristic solution method for the CTP as follows.
First, we find a feasible solution to the CTP that ensures each $w \in W$ is covered once.
Using this solution, we construct a second instance to the CTP.
The second instance is equivalent to the first \emph{except} that (1) the new set of required drop box locations includes all locations selected in the first solution, (2) the locations selected in the first solution are removed from all covering sets $N_w$, and (3) any $w \in W$ that was covered by at least two locations selected in the first solution is removed from $W$. 
A solution to the second CTP instance is guaranteed to be feasible for the DBLP when $r = 0$ and $q = 2$.
A proof of this statement is provided in the Supplementary materials \ref{appx:proofs}.
If $q$ takes a value greater than $2$, it is fairly trivial to extend the process used when $q = 2$.

Once this initial solution corresponding to $r = 0$ is found, we wish to find solutions that are feasible for a larger $r$.
These solutions can be found as follows. 
Initialize $r = 0$.
We iterate and increase the value of $r$ by some predetermined, sufficiently small $\varepsilon$ in each iteration.
Within an iteration, start with the solution identified by the previous iteration.
Identify all pairs of drop box locations such that one is already included in the tour (not belonging to the set $T$) and the other is not.
These pairs represent the locations that can be \emph{swapped} (i.e., remove one from the current solution and add the other). 
We allow the pairs to also represent adding a location to the tour without removing another, or removing a location without adding another.
The latter may be advantageous in cases where a drop box location was added in a previous iteration that makes previously included locations redundant and unnecessary.
There are $\mathcal{O}(|N|^2)$ possible pairs.
We let a pair be \emph{feasible} if the drop box locations obtained after the swap
satisfy constraint set \eqref{model:obj2} according to the current value of $r$ and satisfy the multiple coverage defined by constraint set \eqref{model:basecoverage}.
It is straightforward to check the feasibility of each pair.
Note that we do not consider any pair that results in both an increased cost and lower minimal access function value, since it would directly lead to a dominated solution.
Among all feasible pairs, determine the angle between the incumbent solution and the prospective solution, similar to what was done in \cite{current_median_1994}.
Mathematical details can be found in the Supplementary materials \ref{appx:heur}.
Select the pair that leads to the smallest angle; this incentivizes finding solutions with a lower cost.
Then update the tour based on this pair in a cost minimizing way (e.g., minimum cost removal/insertion or other techniques used for the TSP \cite{gendreau_new_1992}).
Continue until all potential drop box locations have been selected in the solution. This is guaranteed to occur after a finite amount of time since $r$ strictly increases by a fixed amount each iteration and is bounded above by one.
Among all identified solutions, disregard the dominated solutions and return the rest.
It can be checked during each iteration whether each new solution is dominated by a previous solution or if the new solution dominates a previous solution.
When a polynomial time heuristic is used for the TSP and CTP, this heuristic also runs in polynomial time.

\section{Case Study}\label{sec:casestudy}
\begin{figure}[hbt!]
    \centering
\subfloat[Milwaukee block group locations ($W$)\label{fig:estimatedpop}]{%
\includegraphics[width = 0.25\textwidth]{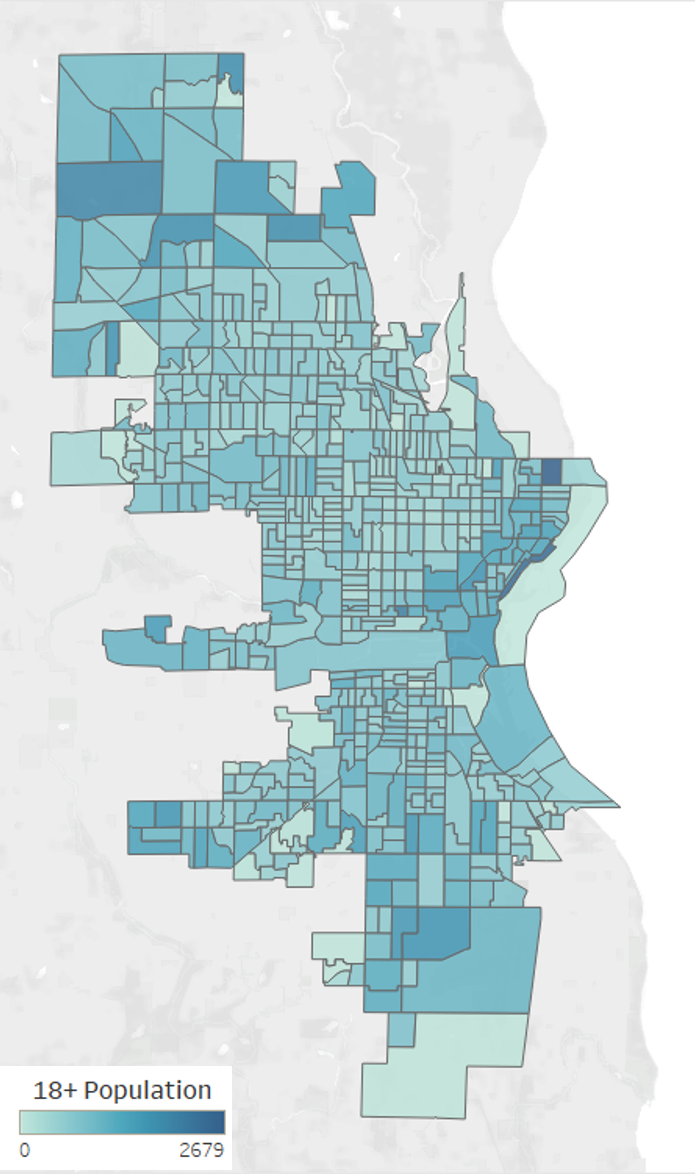}
} %
  \hfill
\subfloat[Drop box locations in 2020 \label{fig:existing_locations}]{%
\includegraphics[width = 0.25\textwidth]{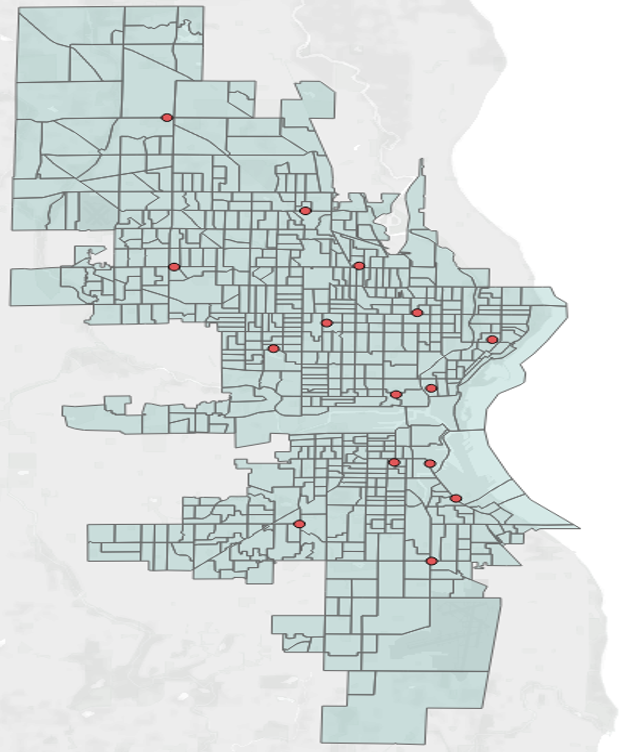}
} %
  \hfill
\subfloat[Potential drop box locations ($N$)\label{fig:potential_locations} ]{%
\includegraphics[width = 0.25\textwidth]{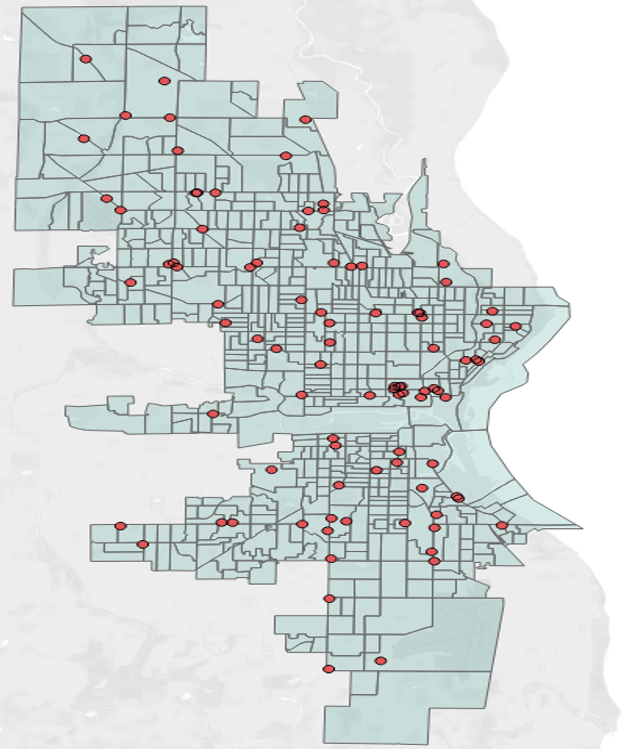}
}%
    \caption{(a) The census block groups of Milwaukee with a darker color indicating a higher population 18 years of age or older.
    (b) The drop box locations (red) used in 2020. (c) Potential drop box locations (red) within the City of Milwaukee used for this case study.}
\end{figure}

We construct a case study of Milwaukee, Wisconsin, United States to demonstrate the value of the DBLP and investigate the implications of optimal drop box infrastructure design.
The City of Milwaukee is the most populous municipality in the state of Wisconsin and had approximately 315,483 registered voters prior to the 2020 General election \cite{city_of_milwaukee_open_data_portal_2020_nodate}. We let $W$ be defined by the census block groups of Milwaukee, WI \cite{milwaukee_county_census_2018},
which are comprised of individuals located near each other who are typically of similar socioeconomic backgrounds. 
Figure \ref{fig:estimatedpop} illustrates the different block group locations in Milwaukee and the estimated number of individuals of age 18 or older in each \cite{united_states_census_bureau_race_2020}.


During the 2020 elections, 15 drop boxes were placed throughout Milwaukee \cite{milwaukee_election_commission_absentee_2020}, illustrated in Figure \ref{fig:existing_locations}.
We use the DBLP to identify drop box locations assuming that these 15 were not already added to the voting system.
This allows us to compare the DBLP to the decisions actually made by election officials during 2020.
We let the potential drop box locations, $N$, be the locations of courthouses (4), election administrative buildings (2), fire stations (30), libraries (14), police stations (7), CVS pharmacies (7), and Walgreens pharmacies (29).
Figure \ref{fig:potential_locations} illustrates the locations of the 93 potential drop box locations.
We assume that the collection tour begins and ends at the Milwaukee City Hall, $s$.
We do not require a drop box be located at any location other than City Hall, with $T = \{s\}$.
The fixed cost 
of locating a drop box at court houses, fire stations, police stations, and City Hall is set at \$6,000 to reflect the cost of a drop box without the need of additional security measures. 
The fixed cost of locating a drop box at all other locations is set at \$10,000 to reflect the cost of both a drop box and a security system.

According to the Milwaukee Election Commission, ballots were, at a minimum, collected daily by staff during the 2020 General election \cite{milwaukee_election_commission_absentee_2020}. This equates to approximately 50 times during the election. 
Based on this value, we assume that ballots are collected 50 times per year on average\footnote{In reality, the number of collections depends on the election year. Also, the frequency of ballot collection may vary depending on model solutions, but this value is set to normalize the operational cost to the fixed cost of the drop boxes. }
over the life of the drop boxes, which we assume to be 15 years.
We further assume that each member of the bipartisan collection team has an opportunity cost of $\$40$ per hour.
This may not reflect the actual pay rate of poll workers or staff; rather, it is  meant to represent the opportunity cost of other tasks not completed during that time.
For example, staff could otherwise participate in additional security training, review compliance of submitted ballots, or conduct marketing to increase voter turnout. 
The cost of traveling between two drop boxes is determined using this pay rate and the estimated time needed to drive between the two locations, which is obtained from Bing Maps.
We include the cost of gas and vehicle wear using the current federal mileage reimbursement of $\$0.56$ per mile.
The estimated mileage is calculated assuming an average travel speed of 30 mph. 
Lastly, we assume the collection costs increase by two percent each year. 

The covering set of each location, $N_w$, is constructed to include the locations that satisfy at least two of the following: the time to walk to the drop box is no more than 15 minutes; the time to drive to the drop box is no more than 15 minutes; the time to use public transit (i.e., city bus) to the drop box is no more than 30 minutes; or the road distance to the drop box is no more than 4 miles (e.g., reachable by bike or ride-share).
By ensuring at least two conditions are are met, there must be multiple transportation modalities that can be used to travel to a drop box in $N_w$.
Individuals without access to a private vehicle are thereby guaranteed to be able to reach a covering drop box using another mode of transportation. 
We estimate the location of each $w \in W$ using the block group centriod.
Throughout the case study, we let $q = 2$, unless otherwise noted, so that model solutions mitigate the risk associated with the destruction of a drop box.

Lastly, the parameters  $v_{w}^0,\ v_w^1,\ a_{nw}$ for the access function are instantiated as follows. 
Ideally, these parameters would be determined using a multinomial logistic regression based on surveys, distance to voting locations, and available transportation methods.
Due to a lack of available data, we introduce a function that serves as a proxy.
Our function combines historical voter turnout, transit durations obtained from Bing Maps, vehicle availability of individuals living in each block group \cite{united_states_census_bureau_means_2019}, and the work locations of individuals residing each block group \cite{united_states_census_bureau_work_nodate}.
We let $v^1_{w}$ equal the estimated voter turnout (between 0 and 1) among registered voters in each block group during the 2016 General election, and we let $v^0_{w} = 1 - v^1_{w}$.
The values reflect the idea that voter turnout is higher where the in-person voting system is more accessible \cite{cantoni_precinct_2020}.
To describe the values of $a_{nw}$, we introduce some notation.
Let $d_{n, w}^\text{walk}$, $d_{n, w}^\text{transit}$, $d_{n, w}^\text{drive}$ be the walking, transit, and driving durations (minutes), respectively, to potential location $n\in N$ for population $w \in W$ obtained from Bing Maps. 
Let $\lambda_{w}^\text{vehicle}$ be the proportion of individuals in $w$ that have access to a vehicle according to the U.S. Census  \cite{united_states_census_bureau_means_2019}.
Let $d_{n, w}^\text{other}$ be the estimated duration to the potential drop box location $n \in N$ from population $w \in W$ using some other form of transportation (e.g., bike or rideshare); the duration is calculated using the road distance obtained from Bing Maps and an assumed speed of 15 miles per hour. 
Let $Q$ represent the set of work locations in Milwaukee, WI \cite{united_states_census_bureau_work_nodate}.
Let $d_{n,q}^\text{work}$ be the walking duration from work location $q \in Q$ to the potential drop box location $n \in N$. 
Let $\lambda_{w,q}^\text{work}$ be the proportion of individuals from $w$ that work in location $q$  according to the U.S. Census \cite{united_states_census_bureau_work_nodate}.
Then, the value of $a_{nw}$ for a drop box location $n$ and population $w \in W$ is computed as:
$$a_{nw} =   \frac{0.04}{v^1_{w}} \Big ( \frac{1}{{(d_{n, w}^\text{walk}})^2} + \frac{1}{({d_{n, w}^\text{transit}})^2} +\frac{\lambda_{w}^\text{vehicle}}{({d_{n, w}^\text{drive}})^2}+\frac{1}{({d_{n, w}^\text{other}})^2} +\sum_{q \in Q}\frac{\lambda_{w,q}^\text{work}}{({d_{n,q}^\text{work}})^2} \Big )$$
This formula accounts for the benefit of multiple modes of transportation to a location and assumes that drop boxes near  individuals are much more accessible than drop boxes far away.
We include the term $(v^1_{w})^{-1}$ to account for the added intangible benefit of drop boxes located near a population with historically low voter turnout, such as increased publicity and visual reminders to cast a ballot \cite{collingwood_drop_2018}.
We square the duration of each transportation mode to model a non-linear relationship between duration and accessibility (a similar approach was employed in \cite{murata_making_2013}).
As a result, $a_{nw}$ is highest when the drop box is a short duration from the voter using each mode of transportation.
We scale each $a_{nw}$ by $0.04$ to produce values that align with findings from  empirical research \cite{mcguire_does_2020}.
Additional explanation and justification of this function is provided in the Supplementary materials \ref{appx:vnw}.

\subsection{The DBLP and Rules-of-thumb}\label{sec:naive}
In the absence of tools to support election planning, election administrators may use rules-of-thumb to select drop box locations.
In this section, we 
demonstrate that the DBLP is able to identify drop box locations that outperform rules-of-thumb across multiple criteria.
The findings support the use of the DBLP during election planning.
Table \ref{T:solutions} presents the details of DBLP solutions for different values of $r$ obtained using the Gurobi 9.1 MIP solver.
Computational studies were run using 64 bit Python 3.7.7 on an Intel\textsuperscript{\textregistered} Core\textsuperscript{TM} i5-7500 CPU with 16 GB of RAM. 
Each optimal solution was identified in less than 3600 seconds.
We refer to the solutions identified by the DBLP as `DBLP $k$' where $k$ refers to the solution ID in Table \ref{T:solutions}.

\begin{table}[!hbt]
\centering 
\def\~{\hphantom{0}}
\caption{Properties of solutions obtained by solving the DBLP with different values of $r$. }
    \label{T:solutions}
          \begin{tabular*}{\columnwidth}{@{}c@{\extracolsep{\fill}}c@{\extracolsep{\fill}}c@{\extracolsep{\fill}}c@{\extracolsep{\fill}}c@{\extracolsep{\fill}}c}
        \toprule
          \shortstack{Solution\\ID} & \shortstack{Minimum Access\\Function Value} & \shortstack{Tour Cost\\ (\$/yr)} & \shortstack{Fixed Cost\\(\$/yr)} & \shortstack{Operational cost\\(\$/yr)} & \shortstack{Number of\\ Drop Boxes} \\
         \midrule
0\lefteqn{^*}&0.573&17,141&6,800&10,341&15\\
1&0.582&17,313&6,800&10,513&15\\
2&0.593&17,813&7,333&10,479&15\\
3\lefteqn{^*}&0.601&18,538&8,267&10,271&16\\
4&0.612&19,701&8,800&10,901&18\\
5&0.620&21,697&10,267&11,430&21\\
6\lefteqn{^*}&0.629&23,599&12,133&11,466&23\\
7&0.637&26,877&14,133&12,743&28\\
8&0.645&31,701&17,200&14,501&33\\
9&0.653&37,318&20,933&16,385&39\\
10&0.661&47,023&27,200&19,823&52\\
11 & 0.668 & 80,574 & 50,800 & 29,774 & 93\\
\bottomrule
\multicolumn{3}{l}{$^*$  \footnotesize  Illustrated in Figure \ref{fig:illustrate_tours}}
\end{tabular*}
\end{table}

During 2020, the Milwaukee Election Commission located drop boxes at the City Hall, the Election Commission warehouse, and 13 neighborhood-based public library branches \cite{milwaukee_election_commission_absentee_2020}.
We begin by comparing these locations to those identified in DBLP 2, which also represents 15 drop boxes.
Table \ref{T:actualcompare} provides the values of multiple criteria for each drop box system.
These criteria provide insight into the performance of each drop box system with respect to cost, access, and risk.

%
\begin{table}[hbt!]
\centering 
\def\~{\hphantom{0}}
    \caption{Comparison of the actual 2020 drop box system to a drop box system design identified by the DBLP across multiple criteria.}
    \label{T:actualcompare}
    \begin{tabular*}{\columnwidth}{@{}p{0.7\columnwidth}|@{\extracolsep{\fill}}c@{\extracolsep{\fill}}c}
        \toprule
          Criteria  & \shortstack{2020} & \shortstack{DBLP 2}\\
         \midrule
Number of Drop Boxes&15&15\\
Fixed Cost (\$/year)&9,733&7,333\\
Operational cost (\$/year)&10,566&10,479\\
Total Cost (\$/year)\symfootnote{}&20,300&17,813\\
Fraction of voters covered by 1 drop box (population weighted)&0.995&1.000\\
Fraction of voters covered by 2 drop boxes (population weighted)\symfootnote{} &0.889&1.000\\
Fraction of voters without access to a car covered by at least two drop boxes by non-driving transit (population weighted)\symfootnote{}&0.941&1.000\\
Minimum access function value\symfootnote{} &0.560&0.593\\
Average access function value (population weighted)&0.776&0.772\\
Maximum road distance to closest drop box&7.634&6.311\\
Maximum road distance to third closest drop box&10.55&9.978\\
Average road distance to closest drop box (population weighted)&1.601&1.679\\
Average road distance to closest 3 drop boxes (population weighted)&2.723&2.486\\
\bottomrule
\multicolumn{3}{l}{ \setcounter{mpfootnote}{0} \footnotesize  \symfootnote{} Objective of the DBLP.} \\
\multicolumn{3}{l}{\footnotesize  \symfootnote{} Required by constraint set \eqref{model:basecoverage} of the DBLP.} \\
\multicolumn{3}{l}{\footnotesize  \symfootnote{} Required by constraint set \eqref{model:basecoverage} given our method of instantiating $N_w$ for each $w \in W$.} \\
\multicolumn{3}{l}{\footnotesize  \symfootnote{} Modeled using constraint set \eqref{model:obj2} in DBLP.}
\setcounter{mpfootnote}{0} 
    \end{tabular*}

\end{table}

The results in Table \ref{T:actualcompare} suggests that the DBLP is able to identify drop box locations that perform better across multiple criteria compared to the rule-of-thumb used by election administrators in Milwaukee during the 2020 election.
We find that with the same number of drop boxes, the DBLP is able to identify drop box locations that result in a lower fixed cost, operational cost, and total cost.
Despite having a lower cost, all voters are covered by at least two drop boxes with DBLP 2, while the 2020 policy only covers 88.9\% of voters twice. 
This gap also exists when voters do not have access to a vehicle (1.00 vs. 0.941).
This means that DBLP 2 admits a higher level of equity of access to the drop box system while mitigating the risk associated with the possible destruction of a drop box.
Moreover, the minimum access function value is higher (0.593 vs. 0.560) for DBLP 2.
With a strict interpretation of the access function, the block group with the lowest turnout is predicted to have a turnout that is 3.3\% higher if the DBLP 2 was implemented rather than the actual locations.
We find that the average access function value is lower for DBLP 2 than the actual implementation; however, the difference is small (0.772 vs. 0.776). 


%

In different situations, other rules-of-thumb may be used by election administrators. 
We compare the DBLP solutions to six other rules-of-thumb that could have potentially been used instead:
\begin{enumerate}[label=Policy \arabic*, leftmargin=2.3cm,topsep=0pt,noitemsep]
    \item Locate drop boxes at the election administrative buildings (2).
    \item Locate drop boxes at the election administrative buildings (2) and police stations (7).
    \item Locate drop boxes at the election administrative buildings (2) and libraries (14).
    \item Locate drop boxes at the election administrative buildings (2), police stations (7), and libraries (14).
    \item Locate drop boxes at the election administrative buildings (2) and fire stations (30).
    \item Locate drop boxes at the election administrative buildings (2), police stations (7), libraries (14), and fire stations (30).
\end{enumerate}
We choose to assess these policies, since they represent the placement of drop boxes at buildings that should be well-distributed throughout the city.
They are not intended to represent a comprehensive list of possible policies.
Table \ref{T:h+esolutions} provides the values of multiple criteria for each rule-of-thumb and DBLP solutions with a similar number of drop boxes.
The results presented in Table \ref{T:h+esolutions} mirror the findings reported in Table \ref{T:actualcompare}; the DBLP identifies drop box locations that are consistently better across multiple criteria than rules-of-thumb with a similar number of drop boxes.
Moreover, most rules-of-thumb are not feasible for the DBLP. Policy 6, which locates 53 drop boxes, is the only rule-of-thumb policy that guarantees that each $w\in W$ is covered by $q=2$ drop box locations.
Meanwhile, the DBLP can find feasible solutions with as few as 15 drop boxes.

\begin{sidewaystable}
\def\arraystretch{0.9}
    \caption{\centering A comparison of rule-of-thumb and DBLP policies across multiple criteria.}
    \label{T:h+esolutions}
    \begin{tabular*}{\columnwidth}{@{}p{6cm}|@{\extracolsep{\fill}}c@{\extracolsep{\fill}}c @{\extracolsep{\fill}}c@{\extracolsep{\fill}}c@{\extracolsep{\fill}}c@{\extracolsep{\fill}}c |@{\extracolsep{\fill}}c@{\extracolsep{\fill}}c@{\extracolsep{\fill}}c@{\extracolsep{\fill}}c}
        \toprule
          Criteria  & \shortstack{Policy\\1} & \shortstack{Policy\\2} & \shortstack{Policy\\3} & \shortstack{Policy\\4} & \shortstack{Policy\\5} & \shortstack{Policy\\6} & \shortstack{DBLP\\3}&  \shortstack{DBLP\\6}& \shortstack{DBLP\\8} & \shortstack{DBLP\\10}\\
         \midrule
Number of Drop Boxes&2&9&16&23&32&53&16&23&33&52\\
Fixed Cost (\$/year)&1,067&3,867&10,400&13,200&13,067&25,200&8,267&12,133&17,200&27,200\\
Operational cost (\$/year)&1,535&7,233&10,616&11,954&18,328&21,129&10,271&11,466&14,501&19,823\\
Total Cost (\$/year)\symfootnote{Objective of the DBLP.}&2,602&11,100&21,016&25,154&31,395&46,329&18,538&23,599&31,701&47,023\\
Fraction of voters covered by 1 drop box (population weighted)&0.596&0.972&0.995&0.995&1.000&1.000&1.000&1.000&1.000&1.000\\
Fraction of voters covered by 2 drop boxes (population weighted)\symfootnote{Required by constraint set \eqref{model:basecoverage} of the DBLP.}&0.362&0.810&0.924&0.973&0.997&1.000&1.000&1.000&1.000&1.000\\
Fraction of voters without access to a car covered by at least two drop boxes by non-driving transit (population weighted)\symfootnote{Required by constraint set \eqref{model:basecoverage} given our method of instantiating $N_w$ for each $w \in W$.}&0.467&0.920&0.963&0.981&1.000&1.000&1.000&1.000&1.000&1.000\\
Minimum access function value\symfootnote{Modeled using constraint set \eqref{model:obj2} in DBLP.}&0.542&0.558&0.568&0.582&0.591&0.623&0.601&0.629&0.645&0.661\\
Average access function value (population weighted)&0.760&0.770&0.777&0.786&0.790&0.809&0.773&0.781&0.788&0.806\\
Maximum road distance to closest drop box&19.269&9.978&7.634&7.634&5.568&5.568&6.311&6.311&6.311&4.865\\
Maximum road distance to third closest drop box&20.123&14.199&10.554&9.978&8.567&8.249&9.851&8.653&8.249&8.249\\
Average road distance to closest drop box (population weighted)&5.829&2.130&1.550&1.469&1.062&0.917&1.689&1.517&1.386&1.105\\
Average road distance to closest 3 drop boxes (population weighted)&6.687&3.348&2.578&2.122&1.774&1.399&2.436&2.052&1.837&1.534\\
\bottomrule
    \end{tabular*}
\end{sidewaystable}

\subsection{Drop Box Trade-offs}\label{sec:optresults}

\newcommand\drawredsquare{%
\begin{tikzpicture}    
    \tikz\draw[red] (0,0) rectangle (0.5ex,0.5ex);  
  \end{tikzpicture}%
}   
In this section, we further investigate DBLP solutions and explore the trade-offs between criteria within the drop box voting system.
We begin by discussing the trade-off between cost and equity of access to all voting pathways (i.e., the minimum access function value).
%
The solutions in Table \ref{T:solutions} suggest there is a substantial trade-off between the cost of the drop box system and the minimum access function value.
However, the marginal increase in cost to achieve an increase in the minimum access function value is not constant.
From DBLP solution 0 to DBLP solution 1 the average cost of a $0.01$ increase of the minimum access function value is $\$186.84$ per year.
From solutions 3 to 4 the average cost of a $0.01$ increase of the minimum access function value is $\$$1,062.25 per year.
From solutions 10 to 11 the average cost of a $0.01$ increase of the minimum access function value is $\$$50,738 per year.
%
%
This highlights the importance of considering the access function within the DBLP. %
When a low cost solution is desirable, an appropriate value for $r$ allows the DBLP to identify drop box locations that admit a much larger minimum access function value for a relatively low increase in cost (e.g., solutions 1-4).
When drop boxes that admit a large minimum access function value are desirable, it is critical to appropriately set $r$, since a small change in $r$ can lead to solutions of substantially different cost (e.g., solutions 8-10).

We next consider the trade-off between equitable access to the drop box system and equitable access to all voting pathways.
Figure \ref{fig:RHS} plots the cost and minimum access function value of multiple optimal solutions when $q$ is zero (\ref{fig:frontier_zero}), one (\ref{fig:frontier_single}), or two (\ref{fig:frontier_double}) with the latter corresponding to the solutions presented Table \ref{T:solutions}. 
When $q = 0$, the DBLP is able to identify drop box locations that substantially increase the minimum access function value for a relatively small cost.
This suggests that there are cost-effective, equitable drop box locations, even when election officials cannot afford to cover each voter with one or two drop boxes.
In general, equitable access to the drop box system and equitable access to all voting pathways are aligned so that access is improved. 
However,
the difference between the curves corresponding to $q = 0$ (\ref{fig:frontier_zero}) and $q = 1$ (\ref{fig:frontier_single}) represents the cost of ensuring equitable access to the drop box system.
In some cases, this cost can be substantial ($\sim \$6,767$ per year).
This demonstrates the trade-off between selecting drop boxes that ensure all voters have access to the drop box system or using the drop boxes to increase the access function value by ``filling in the gaps'' of the in-person voting system.
We also find a substantial difference between the curves corresponding to $q = 1$ (\ref{fig:frontier_single}) and $q = 2$ (\ref{fig:frontier_double}), particularly when a low cost solution is desired.
This suggests that the cost of mitigating risks associated with the destruction of drop boxes through infrastructure design is relatively high and may not be cost effective.
Instead, it may be more cost effective to respond  to an adverse event after it occurs, since the likelihood of this risk occurring is low \cite{scala_evaluating_2021}.


\begin{minipage}{\textwidth}
  \begin{minipage}[t]{0.48\textwidth}
\begin{tikzpicture}
\begin{axis}[
      width=0.8\textwidth,
      height=\textwidth,
      axis lines=middle,
      legend columns = 2,
      scaled y ticks=false,
        legend style={at={(0.05,0.98)},anchor=north west},
        legend cell align={left},
        x label style={at={(axis description cs:0.5,-0.05)},anchor=north},
        y label style={at={(axis description cs:-0.3,.5)},rotate=90,anchor=south},
        ylabel={Cost (\$/year)},
        yticklabel style={
        /pgf/number format/fixed,
        /pgf/number format/precision=5
        },
        xlabel={\shortstack{Minimum access function value}},
        enlargelimits = true,
        legend columns = 1,
        ymin = 0,
        xmin = 0.54,
          ] 
        \addlegendimage{empty legend} \addlegendentry{\hspace{-.6cm}\textbf{q}};
      \addplot[mark=triangle,black,select coords between index={1}{10}] table [y=C,x=P, col sep=comma]{frontier_mke.csv};\addlegendentry{2} \label{fig:frontier_double}
      \addplot[mark=square,red,select coords between index={1}{12}] table [y=C,x=P, col sep=comma]{frontier_single.csv};\addlegendentry{1} \label{fig:frontier_single}
      \addplot[mark=*,orange,select coords between index={1}{14}] table [y=C,x=P, col sep=comma]{frontier_nocov.csv};\addlegendentry{0} \label{fig:frontier_zero}
      \addplot[only marks, mark size=.7pt] coordinates {
      (0.6635,  53000)
      (0.664,  55000)
      (0.6645,  57000)
    };
      %
      %
 \end{axis}
\end{tikzpicture}
\captionof{figure}{\centering Solutions using different values of $q$.} \label{fig:RHS}
  \end{minipage}
\hfill
  \begin{minipage}[t]{0.48\textwidth}
\begin{tikzpicture}
\begin{axis}[
    width=0.8\textwidth,
    height=\textwidth,
    axis lines=middle,
      legend columns = 1,
      scaled y ticks=false,
        legend style={at={(0.05,0.98)},anchor=north west},
        legend cell align={left},
        x label style={at={(axis description cs:0.5,-0.05)},anchor=north},
        y label style={at={(axis description cs:-0.3,.5)},rotate=90,anchor=south},
        ylabel={Cost (\$/year)},
        yticklabel style={
        /pgf/number format/fixed,
        /pgf/number format/precision=5
        },
        xlabel={Minimum access function value},
        enlargelimits = true,
        legend columns = 1,
        ymin = 0,
          ] 
         \addlegendimage{empty legend} \addlegendentry{\hspace{-.6cm}\textbf{Factor}};
      \addplot[orange,mark = otimes,select coords between index={1}{7}] table [y=C_0.9,x=P_0.9, col sep=comma]{frontier_equality.csv};\addlegendentry{0.9} \label{fig:frontier_0.9}
      \addplot[black,mark = triangle,select coords between index={1}{7}] table [y=C_1.0,x=P_1.0, col sep=comma]{frontier_equality.csv};\addlegendentry{1.0} \label{fig:frontier_1.0}
      \addplot[green,mark = square,select coords between index={1}{8}] table [y=C_1.1,x=P_1.1, col sep=comma]{frontier_equality.csv};\addlegendentry{1.1} \label{fig:frontier_1.1}
      \addplot[blue,mark = *,select coords between index={1}{9}] table [y=C_1.2,x=P_1.2, col sep=comma]{frontier_equality.csv};\addlegendentry{1.2} \label{fig:frontier_1.2}
      \addplot[red,mark = diamond,select coords between index={1}{9}] table [y=C_1.3,x=P_1.3, col sep=comma]{frontier_equality.csv};\addlegendentry{1.3} \label{fig:frontier_1.3}
      \addplot[only marks, mark size=.7pt] coordinates {
      (0.66-0.003,  41500)
      (0.661-0.003,  43000)
      (0.662-0.003,  44500)
    };
 \end{axis}
\end{tikzpicture}
\captionof{figure}{\centering Solutions using different time thresholds for $N_w$ determined by the factor.} \label{fig:frontier_compare}
\end{minipage}
\end{minipage}
\vspace{.1cm}

Rather than changing $q$, we can also relax coverage by defining the coverage sets $N_w$ using a larger time threshold.
When covering sets are defined by a longer time threshold, the drop boxes are allowed to be located further away from the voters while still meeting the coverage constraints defined in constraint set \eqref{model:basecoverage}.
A larger time threshold may increase the inequity of access to the drop box infrastructure within the resulting solutions, since a census block group may be further from both covering drop boxes when compared to other census block groups.
However, the access function continues to evaluate the effect of the drop box locations on voter turnout when all voting pathways are considered.
Figure \ref{fig:frontier_compare} illustrates the cost and minimum access function value of solutions to the DLBP using covering sets $N_w$ obtained using the same procedure as before, except the time threshold are multiplied by factors of 0.9 (\ref{fig:frontier_0.9}), 1.0 (\ref{fig:frontier_1.0}), 1.1 (\ref{fig:frontier_1.1}), 1.2  (\ref{fig:frontier_1.2}), and 1.3 (\ref{fig:frontier_1.3}) when $q = 2$.
A factor of 1.0 corresponds to the covering sets used to obtain the solutions discussed in Table \ref{T:solutions} and Figure \ref{fig:RHS}.
We find the effects of changing the covering sets to be similar to the effect of changing $q$.
Covering sets defined by a larger factor result in solutions with a larger minimum access function value for the same cost. 

\begin{figure}[hbt!]
\subfloat[Solution 0]{%
  \includegraphics[width=0.25\columnwidth]{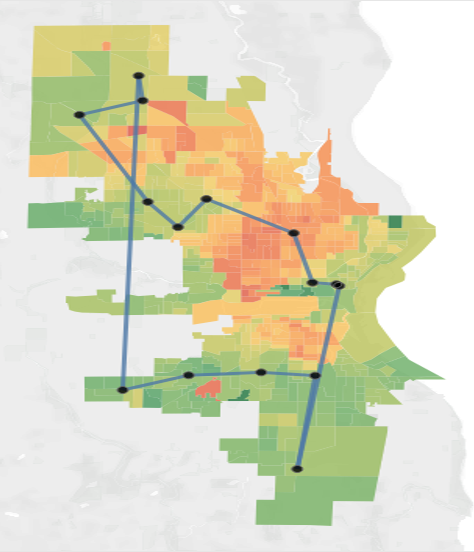}%
}%
\hfill
\subfloat[Solution 3]{%
  \includegraphics[width=0.25\columnwidth]{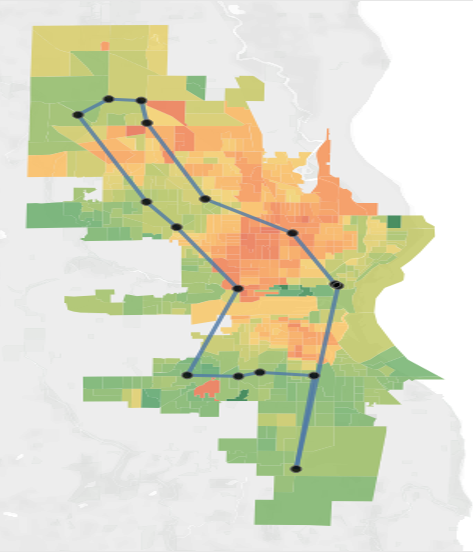}%
}%
\hfill
\subfloat[Solution 6]{%
  \includegraphics[width=0.25\columnwidth]{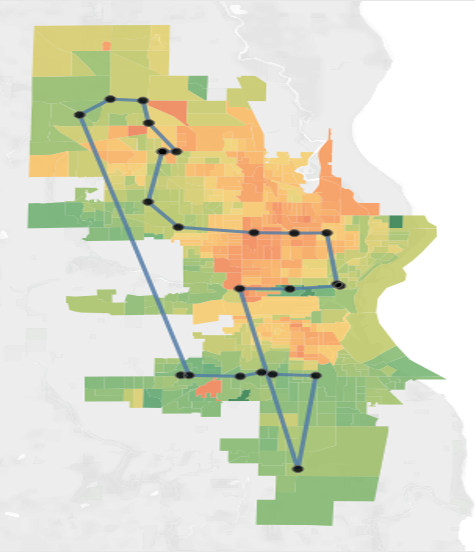}%
}%
\caption{Optimal drop box locations and tour visiting these locations. Color of regions reflect the access function value (red is low, green is high).}\label{fig:illustrate_tours}
\end{figure}

We explore solutions DBLP 0, 3, and 6 of Table \ref{T:solutions} in more detail. Figure \ref{fig:illustrate_tours} illustrates the selected drop box locations and the collection tour visiting these drop boxes overlaid on a map of Milwaukee. 
Each black circle indicates a selected drop box location, and the blue lines describe the order in which the drop boxes are visited on the collection tour (not the actual roads driven).
The color of each block group indicates the access function value of each block group with red reflecting relatively a low value and green reflecting relatively high value.
When cost is minimized (Solution 0), drop boxes are well-spaced in order to cover each block group twice, but relatively few drop boxes are placed to reduce cost.
Solution 0 is notably different than the locations selected during 2020, Figure \ref{fig:existing_locations}, in the northern and southern areas of the city despite locating the same number of drop boxes.
The DBLP selects additional drop box locations in the north and south to ensure equitable access to the drop box system within those region.
When the cost and  minimal access function value are higher (Solutions 3 and 6), more drop boxes are added.
The DBLP selects additional drop box locations in the middle and northern part of the city, which would otherwise have have a relatively low level of access to the voting infrastructure, indicated by the dark red in Figure \ref{fig:illustrate_tours}(a). Additional locations are not selected in the south, since those voters have relatively high access to the multiple voting pathways, indicated by the dark green  in Figure \ref{fig:illustrate_tours}(a).



\subsection{Heuristic Results}
In this section, we investigate the performance of the DBLP heuristic compared to the lazy constraint approach.
The lazy constraints were implemented using the Gurobi 9.1 MIP solver.
Instances to the TSP and CTP created during the heuristic's execution were solved using the following methods.
To solve a CTP instance, we formulate a group Steiner Tree problem instance to determine the nodes to visit in the CTP tour \cite{garg_polylogarithmic_2000} using the approach in \cite{duin_solving_2004} coupled with the
technique introduced in \cite{gubichev_fast_2012}.
Due to the relatively small number of drop box locations, we optimally solve the TSP over selected locations using the Gurobi 9.1 MIP solver. 

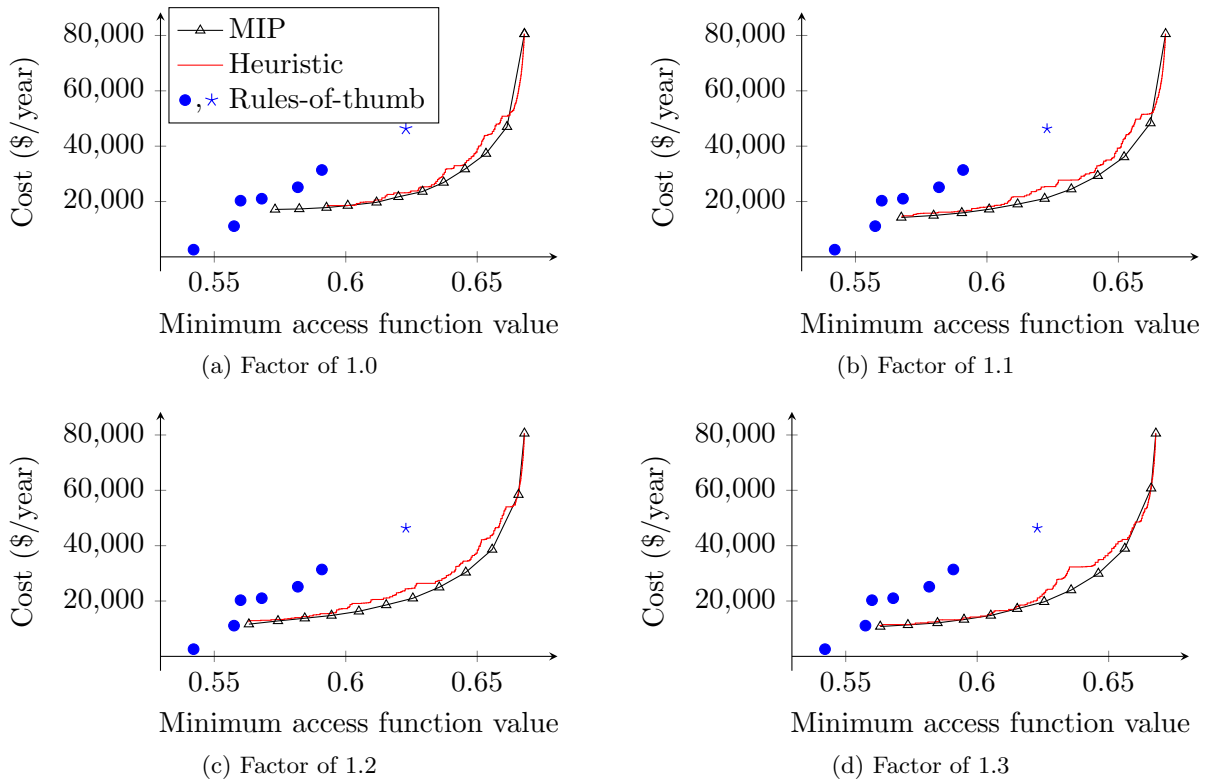
\begin{figure}[htb!]\vspace{0.5cm} 

\subfloat[Factor of 1.0]{%
   \begin{tikzpicture}
\begin{axis}[
	width=0.95*0.45\textwidth,
      height=0.7*0.45\textwidth,
      axis lines=middle,
      legend columns = 2,
      scaled y ticks=false,
        legend style={at={(0.02,1.02)},anchor=north west},
        legend cell align={left},
        x label style={at={(axis description cs:0.5,-0.12)},anchor=north},
        y label style={at={(axis description cs:-0.28,.5)},rotate=90,anchor=south},
        ylabel={Cost (\$/year)},
        yticklabel style={
        /pgf/number format/fixed,
        /pgf/number format/precision=5
        },
        xlabel={Minimum access function value},
        num1/.style={mark size=2pt,mark options={color=blue}, only marks,mark=*},
        num2/.style={mark = star, mark size=2.5pt,mark options={color=blue}, only marks},
        combo legend/.style={
          legend image code/.code={
            \draw [/pgfplots/num1] plot coordinates {(1mm,0cm)};
            \draw plot coordinates {(1.5mm,-3pt)} node {, };
            \draw [/pgfplots/num2] plot coordinates {(4.5mm,0cm)};
          }
        },
        enlargelimits = true,
        legend columns = 1,
          ] 
      \addplot[mark=triangle,black] table [y=C,x=P, col sep=comma]{frontier_mke.csv};\addlegendentry{MIP} \label{fig:frontier_exact}
      \addplot[const plot mark right,red] table [y=C_1.0,x=P_1.0, col sep=comma]{frontier_heuristic_mke.csv};\label{fig:frontier_heuristic}\addlegendentry{Heuristic} 
      %
      \addplot[forget plot, num1] plot coordinates {(0.542119954600902,2601.676216)
        (0.557504137420686,11099.9264514)
        (0.568047601927534,21015.6370626667)
        (0.581764939515681,25154.43166572)
        (0.590925741698136,31394.7444396)
        (0.559995693279841,20299.5882997999)
        }; 
        \addplot[forget plot, num2] plot coordinates 
        {(0.622864340093998,46328.5055448666)
        };

        \addlegendimage{combo legend}
        \addlegendentry{Rules-of-thumb}

 \end{axis}
\end{tikzpicture}
}\label{fig7:a} %
  \hfill
 \subfloat[Factor of 1.1]{%
 \begin{tikzpicture}
\begin{axis}[
	width=0.95*0.45\textwidth,
      height=0.7*0.45\textwidth,
      axis lines=middle,
      legend columns = 2,
      scaled y ticks=false,
        legend style={at={(0.5,-0.38)},anchor=south},
        legend cell align={left},
        x label style={at={(axis description cs:0.5,-0.12)},anchor=north},
        y label style={at={(axis description cs:-0.28,.5)},rotate=90,anchor=south},
        ylabel={Cost (\$/year)},
        yticklabel style={
        /pgf/number format/fixed,
        /pgf/number format/precision=5
        },
        xlabel={Minimum access function value},
        enlargelimits = true,
        legend columns = 3,
          ] 
      \addplot[mark=triangle,black] table [y=C_1.1,x=P_1.1, col sep=comma]{frontier_equality.csv};\addlegendentry{Branch and Bound Method} \label{fig:frontier_exact}
      \addplot[const plot mark right,red] table [y=C_1.1,x=P_1.1, col sep=comma]{frontier_heuristic_mke.csv};\addlegendentry{Heuristic Method} \label{fig:frontier_heuristic}
      \legend{}; 
      %
      \addplot[mark size=2pt,mark options={color=blue}, only marks] plot coordinates {(0.542119954600902,2601.676216)
        (0.557504137420686,11099.9264514)
        (0.568047601927534,21015.6370626667)
        (0.581764939515681,25154.43166572)
        (0.590925741698136,31394.7444396)
        (0.559995693279841,20299.5882997999)
        }; \label{fig:frontier_rot}
         \addplot[mark size=2pt,mark options={color=blue}, only marks,mark = star] plot coordinates {
        (0.622864340093998,46328.5055448666)
        }; 
         \label{fig:frontier_rot_feas}
 \end{axis}
\end{tikzpicture}}\label{fig7:b} %

\subfloat[Factor of 1.2]{%
 \begin{tikzpicture}
\begin{axis}[
	width=0.95*0.45\textwidth,
       height=0.7*0.45\textwidth,
      axis lines=middle,
      legend columns = 2,
      scaled y ticks=false,
        legend style={at={(0.5,-0.38)},anchor=south},
        legend cell align={left},
        x label style={at={(axis description cs:0.5,-0.12)},anchor=north},
        y label style={at={(axis description cs:-0.28,.5)},rotate=90,anchor=south},
        ylabel={Cost (\$/year)},
        yticklabel style={
        /pgf/number format/fixed,
        /pgf/number format/precision=5
        },
        xlabel={Minimum access function value},
        enlargelimits = true,
        legend columns = 3,
          ] 
      \addplot[mark=triangle,black] table [y=C_1.2,x=P_1.2, col sep=comma]{frontier_equality.csv};\addlegendentry{Branch and Bound Method} \label{fig:frontier_exact}
      \addplot[const plot mark right,red] table [y=C_1.2,x=P_1.2, col sep=comma]{frontier_heuristic_mke.csv};\addlegendentry{Heuristic Method} \label{fig:frontier_heuristic}
      \legend{}; 
      %
      \addplot[mark size=2pt,mark options={color=blue}, only marks] plot coordinates {(0.542119954600902,2601.676216)
        (0.557504137420686,11099.9264514)
        (0.568047601927534,21015.6370626667)
        (0.581764939515681,25154.43166572)
        (0.590925741698136,31394.7444396)
        (0.559995693279841,20299.5882997999)
        }; 
        \addplot[mark size=2pt,mark options={color=blue}, only marks,mark = star] plot coordinates {
        (0.622864340093998,46328.5055448666)
        };
 \end{axis}
\end{tikzpicture} 
 }\label{fig7:c} %
  \hfill
  \subfloat[Factor of 1.3]{%
 \begin{tikzpicture}
\begin{axis}[
	width=0.95*0.45\textwidth,
      height=0.7*0.45\textwidth,
      axis lines=middle,
      legend columns = 2,
      scaled y ticks=false,
        legend style={at={(0.5,-0.38)},anchor=south},
        legend cell align={left},
        x label style={at={(axis description cs:0.5,-0.12)},anchor=north},
        y label style={at={(axis description cs:-0.28,.5)},rotate=90,anchor=south},
        ylabel={Cost (\$/year)},
        yticklabel style={
        /pgf/number format/fixed,
        /pgf/number format/precision=5
        },
        xlabel={Minimum access function value},
        enlargelimits = true,
        legend columns = 3,
          ] 
      \addplot[mark=triangle,black] table [y=C_1.3,x=P_1.3, col sep=comma]{frontier_equality.csv}; \label{fig:frontier_exact}\addlegendentry{Branch and Bound Method}
      \addplot[const plot mark right,red] table [y=C_1.3,x=P_1.3, col sep=comma]{frontier_heuristic_mke.csv};\label{fig:frontier_heuristic}\addlegendentry{Heuristic Method} 
      \legend{}; 
      %
      \addplot[mark size=2pt,mark options={color=blue}, only marks] plot coordinates {(0.542119954600902,2601.676216)
        (0.557504137420686,11099.9264514)
        (0.568047601927534,21015.6370626667)
        (0.581764939515681,25154.43166572)
        (0.590925741698136,31394.7444396)
        (0.559995693279841,20299.5882997999)
        };
        \addplot[mark size=2pt,mark options={color=blue}, only marks,mark = star] plot coordinates {
        (0.622864340093998,46328.5055448666)
        };
 \end{axis}
\end{tikzpicture}
 } 
  \caption{\centering Approximate Pareto frontier identified by the heuristic method (\ref{fig:frontier_heuristic}) compared to MIP solutions (\ref{fig:frontier_exact}). Blue dots (\ref{fig:frontier_rot}) represent the cost and minimum access function value of the rules-of-thumb discussed in Section \ref{sec:naive}, 
  with blue stars (\ref{fig:frontier_rot_feas}) indicating feasible solutions for $q=2$.}\label{fig:frontier_heur} 
\end{figure}

Our computational studies suggest that that proposed heuristic method approximates the Pareto frontier between cost and the minimum access function value well and does so quickly.
Figure \ref{fig:frontier_heur} plots the cost and minimum access function of solutions found using the MIP solver (\ref{fig:frontier_exact}) against heuristic solutions  (\ref{fig:frontier_heuristic}) for four instances of the Milwaukee case study, each corresponding to covering sets defined by a different factor (1.0, 1.1, 1.2, or 1.3).
We also plot the cost and minimum access function value of the rules-of-thumb (\ref{fig:frontier_rot}) discussed in Section \ref{sec:naive}; however, recall that these policies may not be feasible for the DBLP.
%
%
We find two primary benefits of the heuristic method.
First, a large number of solutions are identified quickly.
Our experiments showed that the MIP solver can find a new policy every 225 seconds on average when using lazy constraints, while the heuristic method is able to identify between 141-181 policies in 126-220 seconds in total, depending on the factor used to construct the covering sets.
Second, the difference between the cost and minimum access function value of solutions identified by the heuristic and MIP methods is small. 
Moreover, the heuristic is able to identify solutions that are feasible for the DBLP that have a lower cost and higher minimum access function value than rules-of-thumb which may not be feasible for the DBLP.
Election officials can implement the heuristic solutions or use them to determine a range of  appropriate $r$  values and then explore optimal solutions within this range.

We also investigate the heuristic  on randomly generated DBLP instances to further assess its performance. 
The random instances were generated as follows.
The location of voter populations and drop boxes were randomly selected within a 100 by 100 grid. 
We let the time to travel between each location be the Manhattan ($\ell$1) distance.
We randomly select up to a quarter of the drop boxes be required locations ($T$).
Fixed costs were randomly selected between \$5,000 and \$12,000 for each drop box.
Operational costs were computed using the same assumptions described in Section \ref{sec:casestudy}, except costs were scaled by a random value between 0.5 and 1.5.
The time threshold used to construct the covering sets were randomly selected between 15 and 50. 
If a larger threshold was required to ensure each covering set included at least two locations, we select the smallest threshold possible.
The value of $v_w^1$ was randomly selected between 50 and 95 and $v_w^0 = 100-v_w^1$ for each $w \in W$.
The value of $a_{nw}$ was computed using the formula $e^{2.5 - d_{nw}/30}$ where $d_{nw}$ represents the Manhattan distance between the $w$ and $n$ for each $w \in W$ and $n\in N$.

A comparison of solutions generated by the heuristic and MIP solver for nine randomly generated instances is provided in Table \ref{T:heur_results}. 
For each DBLP instance, we use the heuristic to obtain an approximation set for the DBLP. 
For each heuristic solution in the approximation set, we solve the DBLP to optimality using a $r$ equal to the minimum access function value admitted by the solution.
We then compute the cost deviation of the solutions using the average percent difference in the cost of heuristic and optimal solutions corresponding to the same value of $r$.
The results presented in Table  \ref{T:heur_results} are consistent with the findings from the Milwaukee case study.
The heuristic method requires substantially less time than the MIP solver to find solutions corresponding to the same $r$ values.
Moreover, the heuristic method finds solutions with small cost deviations, which indicates that the heuristic finds policies that are similar in cost and access to the optimal solutions. 

\begin{table}[htb!]
    \centering
    \caption{Comparison of the heuristic and MIP implementations for random DBLP instances. The cost deviation is calculated as the average percent difference in total cost for solutions corresponding to the same $r$. Each MIP instance was terminated after 3600 seconds if no optimal was found.}
    \label{T:heur_results}
    \begin{tabular*}{\columnwidth}{@{}c@{\extracolsep{\fill}}c@{\extracolsep{\fill}}c@{\extracolsep{\fill}}c@{\extracolsep{\fill}}c@{\extracolsep{\fill}}c@{\extracolsep{\fill}}c}
        \toprule
          $|W|$ & $|N|$ & $|T|$ & \shortstack{Number of\\Solutions} & \shortstack{Heuristic\\Time (s)} & \shortstack{MIP\lefteqn{^*}\\Time (s)} & \shortstack{Cost\\Deviation} \\
         \midrule
100 & 50& 8 & 59 & 5.22 & 23 & 0.52\%\\
500 & 50 & 11 & 61 & 14.42 & 168 & 0.46\% \\
1000 & 50 & 7 & 53 & 28.05 & 90 & 1.98\% \\
100  & 75 &  4 & 122 & 35.54 &  3,939 & 2.76\%	 \\
500 & 75 &  18 & 109 & 45.91 & 960 & 2.90\%	\\
1000 & 75 & 15 & 93 & 65.33 & 661 & 0.97\%	\\
100  & 100  & 9 & 162 & 84.11 & 23,644 & 4.33\% \\
500 & 100 &  1 & 193 & 173.95 & 349,400$^{**}$ & 8.01\%$^{**}$	\\
1000 & 100 & 12  & 117 & 119.53  & 30,109 & 4.12\%	\\
\bottomrule
    \end{tabular*}
 $^*$  \footnotesize  Ignoring duplicate optimal solutions; 
 $^{**}$ \footnotesize 108 instances terminated after 3600 seconds. 
\end{table}

\section{Conclusion}\label{sec:discussion}
In this paper, we introduce a structured and transparent approach to support the planning of ballot drop box voting systems, particularly for U.S.~voting systems.
We do so by formalizing the drop box location problem (DBLP) to identify a set of optimal drop box locations.
The locations are selected to minimize cost while ensuring voters have access to the drop box system and drop box risks are mitigated. 
Using a real-world case study, we demonstrate that the DBLP identifies drop box locations that consistently outperform rules-of-thumb across multiple criteria.
We also find that the trade-off between criteria is non-trivial and requires careful consideration.

Our research suggests that optimization is an important tool for designing the drop box infrastructure.
Simple guidance for designing drop box systems, such as locating one drop box per 15,000 registered voters, or other rules-of-thumb may be overly-simplistic and can cause election administrators to overlook cost-effective drop box locations that address inequalities within the voting infrastructure.
Strategic drop box locations can reduce the ``cost'' of voting to address inequity within the voting system while ensuring that all voters have equitable access to the drop box system.
Future research can utilize the DBLP to answer additional drop box policy questions and support the drafting of legislation surrounding the use of drop boxes.

We introduce a lazy constraint approach to solve the DBLP to optimality.
Computational experiments show that a single optimal solution to the DBLP can be found relatively quickly using this approach within a state-of-the-art MIP solver for moderately sized problem instances. 
However, the  multiple goals of the DBLP often requires multiple solutions to the DBLP.
We find that this can cause computational times to increase to levels that are unreasonable for practice.
This motivates the need for exceptionally quick solution methods for the DBLP.
We introduce a heuristic for the DBLP, and we demonstrate that the heuristic identifies quality solutions quickly.
Initial attempts at reducing the computational time needed to identify optimal solutions using cutting planes originally introduced for the CTP \cite{gendreau_covering_1997} proved unfruitful.
Future research into the theory of the DBLP is needed to reduce solution times for exact methods.



The DBLP is intended to be a component of a larger suite of tools for supporting election administrators understand, assess, and ultimately design different facets of the voting infrastructure.
Ideally, the DBLP and other operations research tools will eventually be integrated into an online platform designed to support election administrators in all aspect of elections planning.
However, current research is limited.
It largely overlooks the vote-by-mail system and the risks of voting systems.
There is a substantial opportunity for the operations research community to support election planning by appropriately modeling voting systems and voting infrastructure.
Future research is needed to understand the temporal aspects of risk, particularly in the absentee voting process, and determine best practices for mitigating against malicious and non-malicious attacks.
The DBLP and future models can then be incorporated into a comprehensive tool to support election officials in designing the election infrastructure in a way that increases voter turnout.
A key challenge within this space is the need to understand and incorporate models that describe how voters freely select from multiple voting pathways once the infrastructure is set.
Voter choices ultimately determine the cost-effectiveness and performance of the voting system. 



%


\bibliography{refs}

 
%
\appendix
\newpage
\section{Supplementary materials}
\subsection{Proofs}\label{appx:proofs}
\begin{proof}[Proof to Theorem \ref{thm:nphard}]
We reduce the traveling salesman problem (TSP) to the DBLP. Suppose we have an instance of the symmetric TSP defined by the nodes $\bar{N}$, edges $\bar{E}$, and edge costs $\bar{c}$. 
We construct an instance of the DBLP as follows.
Let $N = \bar{N}$, $T = \bar{N}$, $W = \emptyset$, $r = 0$, $q = 0$, $f_j = 0$ for each $j \in \bar{N}$, $E = \bar{E}$, and $c_{ij} = \bar{c}_{ij}$ for $(i,j) \in \bar{E}$.
Then, the DBLP is equivalent to:
\begin{align}
    \underset{x,y}{\min} \ &  z_1 =  \sum_{(i,j) \in \bar{E}} \bar{c}_{ij} x_{ij} \\
    \text{s.t.} \ 
    & \sum_{i \in N : (i,j) \in \bar{E}} x_{ij} = 2 & \forall \ j \in \bar{N}  \\
    & \sum_{i \in S, j \in N\setminus S}    x_{ij} \geq 2  & \forall S \subset \bar{N}, \ 2 \leq |S| \leq |\bar{N}|-2 \\
    & x_{ij} \in \{0,1\} & \forall \ (i, j) \in  \bar{E}
\end{align}
The formulation follows from the fact  that constraint sets \eqref{model:obj2} and \eqref{model:basecoverage} are empty since $W = \emptyset$, and constraint set \eqref{model:existing} requires that $y_{j}$ is to equal one for each $j \in T = \bar{N}$.
This equivalent formulation is an instance of the TSP. Thus, if we can solve the DBLP in polynomial time, then we can solve the TSP in polynomial time. Since the TSP is NP-Hard, so is the DBLP.
\end{proof}

\begin{proof}[Proof to Lemma \ref{prop:covdom}]
Suppose there was a set $N'\subseteq N$, such that $T \subseteq N'$, for which 
$$\frac{v^1_{w}+\sum_{n \in N'} a_{nw}}{v^0_{w}+v^1_{w}+\sum_{n \in N'} a_{nw} } < \frac{v^1_{\hat{w}}+\sum_{n \in N'} a_{n\hat{w}}  }{v^0_{\hat{w}}+v^1_{\hat{w}}+\sum_{n \in N'} a_{n\hat{w}} }$$
Then it must be true that 
$$\frac{v^0_{w} }{v^0_{w}+v^1_{w}+\sum_{n \in N'} a_{nw} } > \frac{v^0_{\hat{w}} }{v^0_{\hat{w}}+v^1_{\hat{w}}+\sum_{n \in N'} a_{n\hat{w}} }$$
This implies that 
$$\frac{v^0_{\hat{w}}+v^1_{\hat{w}}+\sum_{n \in N'} a_{n\hat{w}} }{v^0_{w}+v^1_{w}+\sum_{n \in N'} a_{nw} } > \frac{v^0_{\hat{w}} }{v^0_{w} } \geq 1$$
where the second inequality is true by the assumption of $v_{\hat{w}}^0 \geq v_{w}^0$. This implies
$$v^0_{\hat{w}}+v^1_{\hat{w}}+\sum_{n \in N'} a_{n\hat{w}} > v^0_{w}+v^1_{w}+\sum_{n \in N'} a_{nw} $$
$$ \implies v^0_{\hat{w}}+v^1_{\hat{w}}+\sum_{n \in T} a_{n\hat{w}}+\sum_{n \in N'\setminus T} a_{n\hat{w}} > v^0_{w}+v^1_{w}+\sum_{n \in T} a_{nw}+\sum_{n \in N'\setminus T} a_{nw} $$
$$ \implies (v^0_{\hat{w}}- v^0_{w})+(v^1_{\hat{w}}+\sum_{n \in T} a_{n\hat{w}} -v^1_{w}-\sum_{n \in T} a_{nw})+\sum_{n \in N' \setminus T} (a_{n\hat{w}}-a_{nw}) > 0 $$
However, each  parenthesis term is negative (or zero) by assumption, and the sum can never be greater than zero. This is a contradiction.
\end{proof}


\begin{proof}[Proof to Feasibility of DBLP Solution from Section \ref{sec:heuristic}]
Let  $\hat{x}\in \{0,1\}^{|E|}$ and $\hat{y}\in\{0,1\}^{|N|}$  be the feasible solution to the first CTP instance found.
Let $\hat{T} := T \cup \{n \in N : \hat{y}_n = 1\}$ represent the updated set of required locations used when solving the second CTP instance.
Let $x^*\in \{0,1\}^{|E|}$ and $y^*\in\{0,1\}^{|N|}$ be the feasible solution to the second CTP instance.
Let $N^0 := \{n \in N : y_n^* = 1\}$ denote the drop box locations selected according to $y^*$.
Given a valid solution procedure for the CTP, $x^*$ describes a tour visiting $N^0$.
Thus, the solution must satisfy constraint sets \eqref{model:balance} and \eqref{model:subtourelim} for the DBLP.
Moreover, $T \subseteq \hat{T} \subseteq N^0$ for any feasible solution to the second CTP instance.
Thus, the solution satisfies constraint set \eqref{model:existing} for the DBLP.
What remains to be verified is the satisfaction of constraint set \eqref{model:basecoverage}. By construction:
\begin{align}
\sum_{n \in N_w} y^*_n 
&= |N_w \cap N^0| & \nonumber \\
&= |N_w \cap \hat{T}| + |N_w \cap (N^0\setminus \hat{T})| \nonumber \\
&\geq  1 + |N_w \cap (N^0\setminus \hat{T})| \nonumber \\
&\geq  1 + \sum_{n \in N_w \setminus \hat{T}} y^*_n \nonumber \\
&\geq  1 + 1 \nonumber \\
&\geq  2 \nonumber
\end{align}
The first equality follows from the definition of $N^0$.
The second equality follow from the fact that $\hat{T} \subseteq N^0$. 
The third statement follows from the fact that 
$|N_w \cap \hat{T}| = |N_w \cap \{n \in N : \hat{y}_n = 1\}| \geq 1$.
The fourth statement follows from the definition of $N^0$.
The fifth statement follows from the fact that no location in $\hat{T}$ is a member of the covering sets in the second instance ($N_w \setminus \hat{T}$), and the second CTP must select another location to include within the tour to cover each $w \in W$.

\end{proof}

\subsection{Heuristic Method Psuedocode}\label{appx:heur}
The pseudocode of the DBLP heuristic solution method is presented in Algorithm \ref{alg:heur}.
We assume that the reformulation of the DBLP outlined in Section \ref{sec:obj1reform} is used throughout the heuristic.
We also assume that $q$ takes a value of either one or two; it can easily be extended to cases where $q$ is larger.
When $q = 0$, lines  \ref{alg:CTP1sol}-\ref{alg:CTP2sol} can be replaced so that $x' \in \{0,1\}^{|E|}$ and $y'\in\{0,1\}^{|N|}$ are defined from a (heuristic) solution to the TSP over the set of locations $T$.
The steps of the heuristic are as follows.
The initial solution of drop box locations and the associated collection tour is found in lines \ref{alg:CTP1}-\ref{alg:CTP2sol}.
In lines \ref{alg:CTP1}-\ref{alg:CTP1sol}, a CTP instance is solved based on the DBLP instance.
In lines \ref{alg:CTP2}-\ref{alg:CTP2sol}, a second CTP instance is created and solved.
Solutions to both instances must be feasible for the respective CTP instances, but need not be optimal.
The initial solution for the DBLP is represented by  $N^0$, which represents the selected drop box locations, and 
$\mathcal{C}^0$, which represents the collection tour over the selected locations (lines \ref{alg:DBLP0n}-\ref{alg:DBLP0c}).

In lines \ref{alg:DBLPkstart}-\ref{alg:DBLPkend}, Algorithm \ref{alg:heur} finds DBLP solutions meeting a progressively higher bound $r$ for the minimal access function value.
In line \ref{alg:DBLPkstart}, a set $C$ is initialized to store previously found solutions. This set is used later (lines \ref{alg:checkprev1}-\ref{alg:checkprev2}) determine if Algorithm \ref{alg:heur} has re-found a solution. 
In lines \ref{alg:P}-\ref{alg:Pend}, Algorithm \ref{alg:heur} determines which pairs of drop box locations are feasible by checking constraint sets \eqref{model:obj2} and \eqref{model:basecoverage}.
We do not allow $i = j$, since this would represent no change to the drop box system.
In line \ref{alg:dcost}, we estimate the change in the collection tour cost resulting from the removal of $i$ and insertion of $j$, $\Delta \hat{c}(i,j)$.
This can be estimated using a variety of methods, with the simplest being cheapest cost insertion and shortcut removal.
In line \ref{alg:dr}, we calculate the change in the minimal access function value for pair $(i,j)$, $\Delta r(i,j)$.
In line \ref{alg:value}, we identify the best feasible pair using an angle-based approach similar to that used in \cite{current_median_1994}.
The pairs are assessed based on the angle, using a counter-clockwise orientation, between the vector $\langle -1,0 \rangle$ and the vector ${\langle \Delta r(i,j),\Delta \hat{c}(i,j)  \rangle}$.
The angle, $\theta_{i,j}$ can be calculated using the following formula:
\begin{align}
    \theta(\Delta r, \Delta \hat{c}) =  \begin{cases} 
      2\pi - \cos^{-1}\Big(\frac{-\Delta r}{\sqrt{\Delta r^2+\Delta \hat{c}^2}}\Big) & \Delta \hat{c} \geq 0 \\
      \cos^{-1}\Big(\frac{-\Delta r}{\sqrt{\Delta r^2+\Delta \hat{c}^2}}\Big) & \text{otherwise}
   \end{cases} \label{theta}
\end{align}

  





\noindent Algorithm \ref{alg:heur} selects the pair with the lowest $\theta(\Delta r, \Delta \hat{c})$ in line \ref{alg:bestval}; this leads to DBLP solutions with a lower cost.
In lines \ref{alg:newset}-\ref{alg:newtour}, the incumbent solution is updated based on the selected pair to swap.
In line \ref{alg:checkprev1}-\ref{alg:novelsol}, the minimal access function value, $r$, is updated for the next iteration. 
If Algorithm \ref{alg:heur} has re-found a solution (lines \ref{alg:checkprev1}-\ref{alg:checkprev2}), $r$ is set to be the current minimum access function value.
This avoids future `cycling' where Algorithm \ref{alg:heur}  finds the same solution multiple times.
If Algorithm \ref{alg:heur}  has found a new solution, then  $r$ is updated to be the minimum of $A_w(N^k)$ and  $r+\varepsilon$.
This ensures that the Algorithm \ref{alg:heur} is able to find a feasible solution to the DBLP in the next iteration if one exists (this is a result of $v_{wn} > 0$ for all $n \in N\setminus T$ and $w \in W$).
However, if $\varepsilon$ is sufficiently small, $r$ could be updated by setting $r = r+\varepsilon$.
The value of $\varepsilon$ is sufficiently small when there is guaranteed to be at least one feasible pair to swap in each iteration. 
The following $\varepsilon$ value is guaranteed to be sufficiently small:
$$\varepsilon = \min_{w \in W, n \in N} A_w(N) - A_w(N\setminus n)$$
In lines \ref{alg:end?}-\ref{alg:DBLPkend}, Algorithm \ref{alg:heur}  checks whether to terminate. It terminates when a solution that includes all drop box locations has been found.
Algorithm \ref{alg:heur} returns all non-donminated solutions,  and whether a solution is non-dominated by a new solution can be checked during each iteration.
In an actual implementation, the order of lines within Algorithm \ref{alg:heur} can be optimized to reduce run time (e.g., checking the condition on line 18 before the condition on line 17 may result in shorter run time).

\begin{algorithm}[H]
  \caption{DBLP Heuristic ($\varepsilon$)}\label{alg:heur}
  \begin{algorithmic}[1]
  \Statex \hspace{-0.8cm} \textbf{input} A DBLP instance defined by $(N,T,E,W,N_w,\textbf{c},\textbf{v},\textbf{a},q)$ \label{line:DBLP}
    \Statex \vspace{0.2cm}  \hspace{-0.5cm} \textbf{(CTP)$'$: Find initial CTP solution when $q \geq 1$} 
    \State $(CTP)'$ := instance to the CTP defined by $(N,T,E,W,N_w,\textbf{c})$ \label{alg:CTP1} 
    \State $x' \in \{0,1\}^{|E|},\ y'\in \{0,1\}^{|N|}$ := heuristic solution to $(CTP)'$ \label{alg:CTP1sol}
    \Statex \vspace{0.2cm}  \hspace{-0.5cm} \textbf{(CTP)$''$: Find second CTP solution when $q =2$} 
    \State $W' := W \setminus \{w \in W : |\{n \in N_w : y_n' = 1\}| \geq 2 \}$ \label{alg:CTP2}
    \State $T' := T \cup \{n \in N : y_n' = 1\}$
    \State $N_w' := N_w \setminus T'' \quad \forall \ w \in W'$
    \State $(CTP)''$ := instance to the CTP defined by $(N,T',E,W',N_w',\textbf{c})$
    \State $x'' \in \{0,1\}^{|E|},\ y''\in \{0,1\}^{|N|}$ := heuristic solution to $(CTP)''$ \label{alg:CTP2sol} 
    \Statex \vspace{0.2cm}  \hspace{-0.4cm} \textbf{Initialization}
    \State $N^0 := \{n \in N : y_n = 1\}$ selected locations defined by $y'$ (when $q \leq 1$) or $y''$ (when $q = 2$)  \label{alg:DBLP0n}
    \State $\mathcal{C}^0 :=$ collection tour defined by $x'$ (when $q \leq 1$) or $x''$ (when $q = 2$) \label{alg:DBLP0c}
    \State $C := \{\mathcal{C}^0\} $  a set of previously found solutions \State $r^1 := 0$ \label{alg:DBLPkstart}
    \label{alg:DBLPkprev}
    \algstore{myalg}
    \end{algorithmic}
\end{algorithm}
\begin{algorithm}                     
\begin{algorithmic} [1]                   
\algrestore{myalg}
    \Statex \vspace{0.2cm}  \hspace{-0.4cm} \textbf{Iterative Improvements}
    \For{$k = 1,2,3,\dots, \text{until return} $}
        \State $(i^*,j^*) = \emptyset$
        \State $\theta^* = 2\pi$ 
        \For{$i \in (N^{k-1}\setminus T)$ or $i$ represents no drop box}\label{alg:P}
        \For{$j \in (N\setminus N^{k-1})$ or $j$ represents no drop box ($ i \neq j$)}
        \If{$|N_w \cap (N^{k-1}\cup \{j\}\setminus \{i\})| \geq q \quad \forall \ w \in W$} 
        \If{$r^{k} \leq \min_{w \in W} A_w(N^{k-1}\cup \{j\}\setminus \{i\})$}\label{alg:Pend}
        \State  \label{alg:dcost} \parbox[t]{\dimexpr\textwidth-\leftmargin-\labelsep-\labelwidth}{$\Delta \hat{c} (i,j) = $ estimated change in tour cost from removing $i$ and inserting $j$ }  
        \State \label{alg:dr} \parbox[t]{\dimexpr\textwidth-\leftmargin-\labelsep-\labelwidth}{$\Delta r(i,j) = \min_{w \in W} A_w(N^{k-1}\cup \{j\}\setminus \{i\}) - \min_{w \in W} A_w(N^{k-1})$} 
        \If{$\theta (\Delta r (i,j), \Delta \hat{c} (i,j))$ computed using equation \eqref{theta} $< \theta^*$}\label{alg:value}   
        \State $(i^*,j^*) = (i,j)$ \label{alg:bestval}
        \State $\theta^* =\theta (\Delta r (i,j), \Delta \hat{c} (i,j)) $ 
        \EndIf
        \EndIf
        \EndIf
        \EndFor
        \EndFor
        \State $N^k := N^{k-1} \cup \{i^*\} \setminus \{j^*\}$  \label{alg:newset}
        \State $\mathcal{C}^k := $ heuristic tour over $N^k$   \label{alg:newtour}
        \If{$\mathcal{C}^k \in C$} \label{alg:checkprev1}
            \State $r^{k+1} := \min_{w \in W} A_w(N^k)$ \label{alg:checkprev2}
        \Else
            \State $r^{k+1} := \min \{\min_{w \in W}  A_w(N^k), r+\varepsilon\}$
            \State $C = C \cup \{\mathcal{C}^k\}$ \label{alg:novelsol}
        \EndIf 
        \If{$N^k = N$} \label{alg:end?}
        \State \Return \label{alg:DBLPkend} Identify and return the non-dominated solutions in $C$
        \EndIf
    \EndFor
    \end{algorithmic}
\end{algorithm}

\subsection{Access Function Parameters}\label{appx:vnw}
The access function used throughout this paper is modeled after the conditional/multinomial logit model from discrete choice theory. 
With a strict interpretation of the model, the access function value takes the form
\begin{align}
        A_w(N^*) : = \frac{v^1_{w}+\sum_{n \in N^*} a_{nw} }{v^0_{w}+v^1_{w}+\sum_{n \in N^*} a_{nw} } = \frac{e^{U_w^1}+\sum_{n \in N^*} e^{U_{wn}} }{e^{U_w^0}+e^{U_w^1}+\sum_{n \in N^*} e^{U_{wn}} }  \nonumber 
\end{align}
where $U_w^0$ represents the utility of not voting, $U_w^1$ represents the utility of voting using the non-drop box voting system, and $U_{wn}$ represents the utility of voting by using drop box $n$.
The value of the access function value then represents the probability that an individual chooses to vote using any of the pathways available to them.

According to an economic theory of election participation, potential voters decide whether to vote by comparing the cost to vote and the potential benefits from voting \cite{downs_economic_1957}. 
This idea was later codified as a linear combination of benefits and costs in the form of \cite{riker_theory_1968}
$$\text{Utility} = \text{Benefits} - \text{Costs} $$

\citet{mcguire_does_2020} found that a decrease of one mile to the nearest drop box increases the probability of voting by 0.64 percent.
We use these ideas to justify the method by which we set the value of $a_{nw}$ for all $n \in N$ and $w \in W$, which was presented in Section \ref{sec:casestudy}.
We identify a hypothetical function of the form $a_{nw} \approx e ^{\text{Benefits} - \text{Costs}}$ that aligns with the findings from \citet{mcguire_does_2020}.
We find that $a_{nw} \approx e ^{\text{Benefits} - \text{Costs}} = e^{2.5-D}$ where $D$ is the distance between the voter and the drop box is an appropriate model to validate our method against.
Table \ref{T:utility} is used within our assessment of this hypothetical function.
The second column of Table \ref{T:utility} provides an estimate of the increase in voter turnout for a region when a single drop box a distance of $D$ miles away is added to a voting system that currently has no drop box, assuming ${2.5-D}$ is an appropriate model for the voters' utility.
For example, when a drop box is located a distance of $0.2$ miles from a voter, the expected increase in voter turnout is $2.7\%$.
We assume a 70\% turnout in prior elections for this region.
The values in the second column are calculated as follows
$$\frac{0.7+ e^{2.5-D} }{1+ e^{2.5-D}} - \frac{0.7}{1}$$
where the distance $D$ is given in the first column.
The values in the third column represent the estimated impact of a  one mile decrease to the nearest drop box.
The values are calculated by taking the value in the second column for $D$ and subtracting the value of the second column for a distance that is one mile longer ($1+D$).
For example, the value in the first row is calculated by computing $0.027-0.011$, where the first value corresponds to $D = 0.2$ and the second corresponds to $D = 1.2$.
The average of the values in the third column is 0.0061.
This roughly aligns with the 0.64 percent found by \citet{mcguire_does_2020} as desired.

\begin{table}[H]
    \centering
    \caption{Implications of $e^{2.5-D}$ description of the $a_{nw}$'s assuming 70\% turnout in the voting system without any drop boxes.}
    \label{T:utility}
    \begin{tabular*}{\columnwidth}{@{}c@{\extracolsep{\fill}}c@{\extracolsep{\fill}}c }
        \toprule
          \shortstack{Distance to \\
          Drop  Box\\ 
           (mi), D} & \shortstack{Marginal Increase\\ in Access \\ Function Value} & \shortstack{Benefit of 1 mile \\ Decrease to Drop Box\\ (Resulting in D)}  \\
         \midrule
0.2&0.027&0.017\\
0.4&0.023&0.014\\
0.6&0.019&0.012\\
0.8&0.016&0.010\\
1.0&0.013&0.008\\
1.2&0.011&0.007\\
1.4&0.009&0.005\\
1.6&0.007&0.005\\
1.8&0.006&0.004\\
2.0&0.005&0.003\\
2.2&0.004&0.003\\
2.4&0.003&0.002\\
2.6&0.003&0.002\\
2.8&0.002&0.001\\
3.0&0.002&0.001\\
\bottomrule
    \end{tabular*}
\end{table}

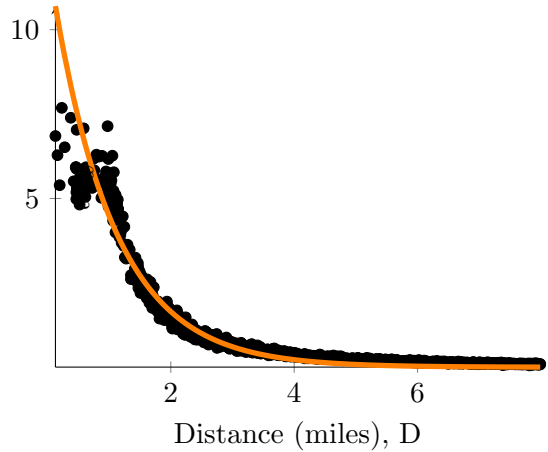
\begin{figure}[ht!]
    \centering
    \begin{tikzpicture}
    \begin{axis}[
            xlabel= {Distance (miles), D}, 
            ylabel=$a_{nw}$, 
            axis lines=left, 
            clip=false,
            clip mode=individual,
            width=0.5\textwidth,
            height=0.4\textwidth,
            y label style={at={(axis description cs:0.1,.5)}},
        ]
                 \addplot[color=black,only marks,mark=*] table [x=D, y=v, col sep=comma] {vnw_data.csv}; \label{mark:sample}
                 \addplot[color=orange, line width=2pt] table [x=D, y=e, col sep=comma] {vnw_data.csv};\label{mark:utility}
        \end{axis}
    \end{tikzpicture}
    \caption{The $a_{nw}$ value and distance between $n$ and $w$ for a sample of 1193 pairs (\ref{mark:sample}) of $n\in N$ and $w \in W$ overlaid with the hypothetical $e^{2.5 - D}$ (\ref{mark:utility}) for which $D <8$.} 
    \label{fig:utilitymodel}
\end{figure}

We now validate the $a_{nw}$ values used in Section \ref{sec:casestudy}.
In Figure \ref{fig:utilitymodel} we plot (\ref{mark:sample}) the $a_{nw}$ value for 1,193 randomly sampled pairs $n \in N$ and $w \in W$ against the distance, $D$, between $n$ and $w$.
Overlaid on these points (\ref{mark:utility}) is the function $e^{2.5-D}$ where the cost is the distance $D$ between $n$ and $w$. 
We find that the proposed method from Section \ref{sec:casestudy} produces $a_{nw}$ values that roughly align with the function $e^{2.5-D}$.
There is variance from the hypothetical line, especially with smaller distances.
This is because we consider additional modes of transit and other factors in the actual calculation of $a_{nw}$. This is desired as it adds more realism. 


\end{document}